\documentclass{article}
\usepackage[leqno,intlimits]{amsmath}
\usepackage{amssymb}
\usepackage{amsthm}
\usepackage{hyperref}
\allowdisplaybreaks[1]

\newtheorem{tm}{Theorem}[section]
\newtheorem{lm}[tm]{Lemma}
\newtheorem{cor}[tm]{Corollary}

\newcommand{\va}{a}  

\newcommand*{\hop}{\bigskip\noindent}
\newcommand{\nn}{\nonumber}
\newcommand*{\Rb}{\mathbb R}
\newcommand*{\Zb}{\mathbb Z}
\newcommand*{\un}{\underline}
\newcommand*{\ba}{\begin{aligned}}
\newcommand*{\ea}{\end{aligned}}
\newcommand*{\be}{\begin{equation}}
\newcommand*{\ee}{\end{equation}}
\newcommand*{\vr}{\varrho}
\newcommand*{\vp}{\varphi}
\newcommand*{\ve}{\varepsilon}
\newcommand*{\om}{\omega}
\newcommand*{\al}{\alpha}
\newcommand*{\la}{\lambda}
\newcommand*{\pt}{\partial}
\newcommand*{\di}{\,\text{\rm d}}

\newcommand*{\Ac}{\mathcal A}
\newcommand*{\Wc}{\mathcal W}
\newcommand*{\Nc}{\mathcal N}
\newcommand*{\Ec}{\mathcal E}
\newcommand*{\Sc}{\mathcal S}
\newcommand*{\wt}{\widetilde}
\newcommand*{\wh}{\widehat}
\newcommand*{\Ev}{\mathbf E}
\newcommand*{\Pv}{\mathbf P}
\newcommand*{\Vv}{{\text{\bf Var}}}
\newcommand*{\Cov}{\text{\bf Cov}}
\newcommand*{\e}[1]{\text{\rm e}^{#1}}
\newcommand*{\fl}[1]{\left\lfloor{#1}\right\rfloor}
\newcommand*{\tfl}[1]{\lfloor{#1}\rfloor}
\newcommand*{\si}{\sigma}
\DeclareMathOperator\Expd{Exp}
\numberwithin{equation}{section}
\begin{document}
\title{Cube root fluctuations for the corner growth model
associated to the
exclusion process}
\author{M.\ Bal\'azs\thanks{University of Wisconsin-Madison\newline
M. Bal\'azs was partially supported by the Hungarian Scientific Research Fund
(OTKA) grant TS49835 and by National Science Foundation Grant DMS-0503650.\newline
T.\ Sepp\"al\"ainen was partially supported by National Science Foundation
grant DMS-0402231.}, E. Cator\thanks{Delft University of Technology}
 and  T.\ Sepp\"al\"ainen$^*$}
\maketitle
\begin{abstract}
We study the last-passage growth model on the planar integer lattice
 with exponential weights.  With boundary conditions
 that represent the equilibrium exclusion process as seen
from a particle right after its jump
we prove that the variance of the last-passage time
in a characteristic direction
is of order $t^{2/3}$.  With more general boundary
conditions that include the rarefaction fan case
we show that the last-passage time fluctuations are
still  of order $t^{1/3}$, and also that
the transversal fluctuations of the maximal path
have order $t^{2/3}$.
We   adapt and then build on   a recent study
of Hammersley's process by Cator and Groeneboom, and also
utilize   the competition interface introduced by Ferrari,
Martin and Pimentel.
  The arguments are entirely probabilistic, and no
use is made of the combinatorics of Young tableaux or methods of
asymptotic analysis.
\end{abstract}

\noindent
{\bf Keywords:} Last-passage, simple exclusion, cube root asymptotics, competition interface, Burke's theorem, rarefaction fan

\hop
{\bf MSC:} 60K35, 82C43

\section{Introduction}

We construct a version of the corner growth model that corresponds
to an equilibrium exclusion process as seen by a typical particle
right after its jump,
and show that along a characteristic direction the variance of the
last-passage time is of order $t^{2/3}$.  This last-passage time
is the maximal sum of  exponential weights along up-right paths
in the first quadrant of the integer plane.
The interior weights have rate 1, while the boundary weights
on the axes have rates
$1-\vr$ and $\vr$ where $0<\vr<1$ is the particle density of the
exclusion process.
By comparison to this equilibrium setting, we also show
fluctuation results with
 similar scaling in the case of the rarefaction fan.

The proof is based on a recent work of Cator and Groeneboom
\cite{cuberoot} where corresponding results are proved for
the planar-increasing-path version of
Ham\-mers\-ley's process.  A key part of that proof is an identity
that relates the variance of the last-passage time to the
point where the maximal path exits the axes. This exit point itself
is related to a second-class particle via a time reversal.
The idea that the current and  the second-class particle
should be connected goes  back to a paper of Ferrari and Fontes
\cite{se} on the diffusive fluctuations
of the current away from the characteristic.  However, despite this
surprising congruence of ideas, article  \cite{cuberoot} and our
work have  no technical relation to the
Ferrari-Fontes work.

The first task of the present paper is to find the
 connection between the variance of the last-passage time and
the  exit point, in the equilibrium corner growth model.  The relation
turns out not as straightforward as for Hammersley's process,
for we also need to include the amount of weight collected on
the axes.  However, once this difference is understood, the arguments
proceed quite similarly to those in  \cite{cuberoot}.

The notion of competition interface recently introduced
by Ferrari, Martin and Pimentel \cite{fermarpim, compint} now appears as
the representative of a second-class particle, and as the time
reversal of the maximal path.  As a by-product of the proof we
establish that the transversal fluctuations of the competition
interface are of the order
$t^{2/3}$ in the equilibrium setting.

In the last section we take full advantage of our
probabilistic approach, and show that for initial conditions obtained
by decreasing the equilibrium  weights on the axes
in an {\sl arbitrary} way,
 the fluctuations
of the last-passage time are still of order $t^{1/3}$.
This includes the situation known as the {\sl rarefaction fan}.
We are also able
to show that in this case the transversal fluctuations of the longest
path are of order $t^{2/3}$. In this more general setting
there is no direct connection between a maximal path and a competition
interface (or trajectory of a second class particle).

Our results for the competition interface, and
our fluctuation results  under the
more general boundary conditions are new.
The variance bound for the equilibrium last-passage time
is also strictly speaking new. However,
 the corresponding
distributional limit has been obtained by Ferrari and Spohn
\cite{ferspohn} with a  proof based on the RSK machinery.
But they lack a suitable tightness property that would give them
 also control of the variance. [Note that
Ferrari and Spohn start by describing a different set of
equilibrium boundary conditions than the ones we consider, but later
in their paper they cover also the kind we  define in
\eqref{eq:bondis} below.]
The methods of our paper can  also be applied to
geometrically distributed weights, with the same outcomes.


In addition to the results themselves,
our main  motivation  is to investigate
 new methods to attack the last-passage model, methods
that do not rely on the RSK correspondence of Young tableaux.
The reason for such a pursuit is that the precise counting
techniques of Young tableaux appear to work only for geometrically
distributed weights, from which one can then take a limit
to  obtain the case of  exponential weights. New techniques are needed
to go beyond the geometric and exponential cases,
although we are not yet in a position to undertake such an advance.


For the  class of totally asymmetric stochastic interacting
systems for which the last-passage approach  works, this point of view
 has been extremely valuable.  In addition to the papers already
mentioned above, we list
 Sepp\"al\"ainen \cite{timogrowth, hkl}, Johansson \cite{1/3}, and
 Pr\"ahofer and Spohn  \cite{spohn}.

{\bf Organization of the paper.}  The main results are discussed
in Section \ref{sc:recults}. Section \ref{sc:part}
describes the relationship of the last-passage model to
particle and deposition models, and can be skipped without loss
of continuity.
The remainder of the paper is for the proofs.  Section \ref{sec:prel}
covers some preliminary matters. This includes a strong form
of Burke's theorem for the last-passage times (Lemma \ref{lm:NE}).
Upper and lower bounds for the equilibrium results are
covered in Sections \ref{sc:ub} and \ref{sc:lb}.
Lastly, fluctuations under more
general boundary conditions are studied in Section \ref{sc:gen-bd}.

{\bf Notation.}  $\Zb_+=\{0,1,2,\dotsc\}$ denotes the set
of nonnegative integers.  The integer part of a real number
is $\tfl{x}=\max\{ n\in\Zb: n\leq x\}$.
$C$ denotes constants whose precise value is immaterial and that do not depend
on the parameter (typically $t$) that grows.
$X\sim \Expd(\vr)$ means that $X$ has the exponential distribution with rate
$\vr$, in other words has density  $f(x)=\vr e^{-\vr x}$ on $\Rb_+$.
  For clarity,
subscripts can be replaced
by arguments in parentheses, as for example in $G_{ij}=G(i,j)$.

\section{Results}
\label{sc:recults}
We start by describing the corner growth model with boundaries
that correspond to a special view of the equilibrium.
Section \ref{sc:part} and Lemma \ref{lm:NE} justify the term equilibrium in this context. Our results for more general
 boundary conditions are in Section \ref{sc:rareres}.

\subsection{Equilibrium results}

  We are given an array
$\{\om_{ij}\}_{i,j\in\Zb_+}$ of nonnegative real numbers.
We will always have $\om_{00}=0$.
The values $\om_{ij}$ with either $i=0$ or $j=0$ are
the boundary values, while $\{\om_{ij}\}_{i,j\geq 1}$ are
the interior values.

\begin{figure}[ht]
\begin{center}
\begin{picture}(220,140)(5,20)
\put(40,30){\vector(1,0){150}}
\put(50,20){\vector(0,1){130}}
\multiput(50,40)(0,20){6}{\line(1,0){140}}
\multiput(60,30)(20,0){7}{\line(0,1){120}}
\put(47,27){\large$\star$}
\multiput(66,26)(20,0){6}{\large$\triangledown$}
\multiput(45.6,47)(0,20){5}{\large$\vartriangle$}
\multiput(67,47)(0,20){5}{
\multiput(0,0)(20,0){6}{\large$\circ$}}
\put(32,23){\tiny$(0,0)$}
\put(185,23){\tiny$i$}
\put(43,145){\tiny$j$}
\put(69,21){\tiny$1$}
\put(89,21){\tiny$2$}
\put(109,21){\tiny$3$}
\put(129,21){\tiny$4$}
\put(149,21){\tiny$5$}
\put(169,21){\tiny$6$}
\put(42,47){\tiny$1$}
\put(42,67){\tiny$2$}
\put(42,87){\tiny$3$}
\put(42,107){\tiny$4$}
\put(42,127){\tiny$5$}
\put(142,63){\tiny$(5,2)$}
\put(142,43){\tiny$(5,1)$}
\put(142,83){\tiny$(5,3)$}
\put(122,63){\tiny$(4,2)$}
\put(162,63){\tiny$(6,2)$}
{\linethickness{2pt}
\put(50,40){\line(1,0){10}}
\put(60,30){\line(0,1){10}}
}
\end{picture}
\end{center}
\caption{The initial situation}\label{fig:iniri}
\end{figure}
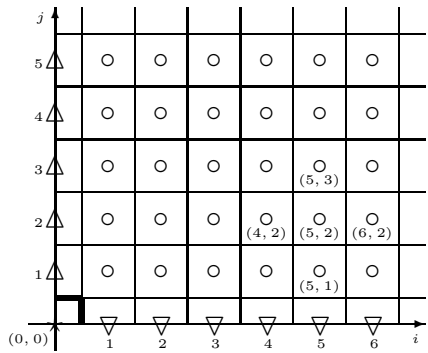

 Figure \ref{fig:iniri} depicts this
initial set-up on the first quadrant $\Zb_+^2$ of the integer plane.
A $\star$ marks $(0,\,0)$,
$\triangledown$'s mark positions $(i,\,0)$, $i\geq1$, $\vartriangle$'s
positions  $(0,\,j)$, $j\geq1$, and interior points $(i,\,j)$, $i,j\geq1$
are marked  with $\circ$'s.  The coordinates of a few
points around $(5,\,2)$ have been labeled.

For a point $(i,\,j)\in\Zb_+^2$, let $\Pi_{ij}$
be the set of directed paths
\be
\pi=\{(0,\,0)=(p_0,\,q_0)\to(p_1,\,q_1)\to\dots\to(p_{i+j},\,q_{i+j})=(i,\,j)\}
\label{eq:pi-1} \ee
with up-right steps
\be
(p_{l+1},\,q_{l+1})-(p_l,\,q_l)= (1,0) \ \text{ or } \ (0,1)
\label{eq:allow}
\ee
along the coordinate directions.
Define the \emph{last passage time} of the point $(i,\,j)$ as
\[
G_{ij}=\max_{\pi\in\Pi_{ij}}\sum_{(p,q)\in\pi}\om_{pq}.
\]
 $G$ satisfies the recurrence
\be
G_{ij}=(G_{\{i-1\}j}\lor G_{i\{j-1\}})+\om_{ij}\qquad(i,\,j\geq0)\label{eq:grec}
\ee
(with formally assuming $G_{\{-1\}j}=G_{i\{-1\}}=0$).
A common interpretation is that this models a  growing cluster
on the
first quadrant that starts from a seed at the origin
(bounded by the thickset line in Figure \ref{fig:iniri}).
The value  $\om_{ij}$
is the time it takes to occupy point $(i,j)$  after its neighbors
to the left and below have become occupied, with the interpretation
that a boundary point needs only one occupied neighbor.
Then $G_{ij}$ is the time when $(i,j)$ becomes occupied, or
 joins the growing cluster.
The occupied region at time  $t\geq0$ is the set
\be
\Ac(t)=\{(i,j)\in\Zb_+^2\,:\, G_{ij}\leq t\}.
\label{eq:def-Ac}
\ee

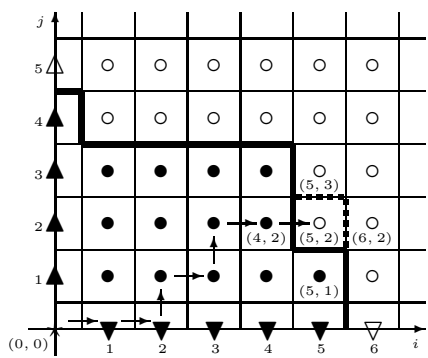
\begin{figure}[ht]
\begin{center}
\begin{picture}(220,140)(5,20)
\put(40,30){\vector(1,0){150}}
\put(50,20){\vector(0,1){130}}
\multiput(50,40)(0,20){6}{\line(1,0){140}}
\multiput(60,30)(20,0){7}{\line(0,1){120}}
\put(47,27){\large$\star$}
\multiput(66,26)(20,0){5}{\large$\blacktriangledown$}
\put(166,26){\large$\triangledown$}
\multiput(45.6,47)(0,20){4}{\large$\blacktriangle$}
\put(45.6,127){\large$\vartriangle$}
\multiput(67,47)(0,20){3}{
\multiput(0,0)(20,0){4}{\large$\bullet$}}
\put(147,47){\large$\bullet$}
\multiput(67,107)(0,20){2}{
\multiput(0,0)(20,0){6}{\large$\circ$}}
\multiput(147,67)(0,20){2}{\large$\circ$}
\multiput(167,47)(0,20){3}{\large$\circ$}
\put(32,23){\tiny$(0,0)$}
\put(185,23){\tiny$i$}
\put(43,145){\tiny$j$}
\put(69,21){\tiny$1$}
\put(89,21){\tiny$2$}
\put(109,21){\tiny$3$}
\put(129,21){\tiny$4$}
\put(149,21){\tiny$5$}
\put(169,21){\tiny$6$}
\put(42,47){\tiny$1$}
\put(42,67){\tiny$2$}
\put(42,87){\tiny$3$}
\put(42,107){\tiny$4$}
\put(42,127){\tiny$5$}
\put(142,63){\tiny$(5,2)$}
\put(142,43){\tiny$(5,1)$}
\put(142,83){\tiny$(5,3)$}
\put(122,63){\tiny$(4,2)$}
\put(162,63){\tiny$(6,2)$}
\multiput(55,33)(60,37){2}{
\multiput(0,0)(20,0){2}{\vector(1,0){11}}}
\multiput(90,35)(20,20){2}{\vector(0,1){10}}
\put(95,50){\vector(1,0){11}}
{\linethickness{2pt}
\put(50,120){\line(1,0){10}}
\put(60,100){\line(0,1){20}}
\put(60,100){\line(1,0){80}}
\put(140,60){\line(0,1){40}}
\put(140,60){\line(1,0){20}}
\put(160,30){\line(0,1){30}}
\multiput(141,80)(4,0){5}{\line(1,0){2}}
\multiput(160,61)(0,4){5}{\line(0,1){2}}
}
\end{picture}
\end{center}
\caption{A possible  later situation}\label{fig:laterri}
\end{figure}

Figure \ref{fig:laterri} shows a possible later situation.
Occupied points are denoted by solidly colored symbols, the
occupied cluster is bounded by the thickset line, and the
arrows mark an admissible path $\pi$ from $(0,0)$ to $(5,2)$.
If $G_{5,2}$ is the smallest among
$G_{0,5}$, $G_{1,4}$, $G_{5,2}$ and $G_{6,0}$, then $(5,2)$
is the next point added to the cluster, as suggested by the
dashed lines around the $(5,2)$ square.

To create a model of random evolution, we pick  a real number
$0<\vr<1$ and  take  the variables $\{\om_{ij}\}$  mutually independent
with the following marginal distributions:
\be
\ba
\om_{00}&=0,\qquad&&\text{where the }\star\text{ is},\\
\om_{i0}&\sim\Expd(1-\vr),\ i\geq1,\qquad&&\text{where the }\triangledown\text{'s are},\\
\om_{0j}&\sim\Expd(\vr),\ j\geq1,\qquad&&\text{where the }\vartriangle\text{'s are},\\
\om_{ij}&\sim\Expd(1),\ i,j\geq1,\qquad&&\text{where the }\circ\text{'s are}.
\ea\label{eq:bondis}
\ee

Ferrari, Pr\"ahofer and Spohn \cite{ferspohn}, \cite{spohn} consider the Bernoulli-equilibrium of simple exclusion, which corresponds to a slightly more complicated boundary distribution than the one described above. However, Ferrari and Spohn \cite{ferspohn} early on turn to the distribution described by \eqref{eq:bondis}, as it is more natural for last-passage. We will greatly exploit the simplicity of \eqref{eq:bondis} in Section \ref{sec:prel}. In fact, \eqref{eq:bondis} is also connected with the stationary exclusion process of particle density $\vr$. To see this point, we need to look at a particle of simple exclusion in a specific manner that we explain below in Section \ref{sec:tasep}.

Once the parameter $\vr$ has been picked we denote the last-passage
time of point $(m,n)$ by $G^\vr_{mn}$.
In order to see interesting behavior we follow
the last-passage time  along the  ray defined by
\be
(m(t),n(t))=(\tfl{(1-\vr)^2t}, \tfl{\vr^2t})
\label{eq:char-mn}
\ee
 as $t\to\infty$.  In Section \ref{sec:char}   we give a
heuristic justification for this choice. It represents the
characteristic speed of the macroscopic equation of the system.   Let us
abbreviate
\[
G^\vr(t)=G^\vr\bigl(\tfl{(1-\vr)^2t}, \tfl{\vr^2t} \bigr).
\]
Once we have proved that all horizontal and vertical increments
of $G$-values are distributed exponentially like the boundary increments,
 we see that
\[ \Ev(G^\vr(t))=\frac{\tfl{(1-\vr)^2t}}{1-\vr}+\frac{\tfl{\vr^2t}}{\vr}.
\]
The first result is the order of the variance.

\begin{tm}\label{tm:G-var-1}
With $0<\vr<1$
and  independent $\{\om_{ij}\}$ distributed as in {\rm \eqref{eq:bondis}},
\[
0<\liminf_{t\to\infty} \frac{\Vv(G^\vr(t))}{t^{2/3}} \leq
\limsup_{t\to\infty} \frac{\Vv(G^\vr(t))}{t^{2/3}} <\infty.
\]
\end{tm}

For given $(m,n)$  there is almost
surely a unique path $\wh\pi$ that maximizes the
passage time to $(m,n)$,
due to the continuity of the distribution of
$\{\om_{ij}\}$.
The {\sl exit point} of $\wh\pi$ is the last boundary point
on the path. If $(p_l,q_l)$ is the exit point for the path
in  \eqref{eq:pi-1}, then either $p_0=p_1=\dotsm=p_l=0$
or  $q_0=q_1=\dotsm=q_l=0$, and $p_k,q_k\geq 1$ for all $k>l$.
To distinguish between exits via the $i$- and $j$-axis, we introduce
a non-zero integer-valued random variable $Z$ such that if
$Z>0$ then the exit point is $(p_{\lvert{Z}\rvert},q_{\lvert{Z}\rvert})=(Z,0)$,
while  if
$Z<0$ then the exit point is $(p_{\lvert{Z}\rvert},q_{\lvert{Z}\rvert})
=(0,-Z)$.  For the sake of convenience we abuse language and call
the variable $Z$ also the ``exit point.''
 $Z^\vr(t)$ denotes the exit point of the maximal path
to the point $(m(t),n(t))$ in \eqref{eq:char-mn}
 with
boundary condition parameter $\vr$.
Transposition $\om_{ij}\mapsto\om_{ji}$ of the array shows
that $Z^\vr(t)$ and $-Z^{1-\vr}(t)$ are equal in distribution.
Along the way to Theorem \ref{tm:G-var-1} we establish
that $Z^\vr(t)$ fluctuates on the scale $t^{2/3}$.

\begin{tm} Given $0<\vr<1$ and independent $\{\om_{ij}\}$ distributed
 as in {\rm \eqref{eq:bondis}}.

{\rm (a)} For $t_0>0$
there exists a finite constant  $C=C(t_0, \vr)$
such that, for all $a>0$ and $t\geq t_0$,
\[
\Pv\{ Z^\vr(t)\geq at^{2/3}\} \leq C a^{-3}.
\]

{\rm (b)} Given $\ve>0$, we can choose a
 $\delta>0$ small enough so that for all large enough  $t$
\[
\Pv\{1\leq  Z^\vr(t)\leq \delta t^{2/3}\}\leq \ve.
\]
\label{tm:Z-1}
\end{tm}

\hop
{\bf Competition interface.}
In \cite{fermarpim, compint} Ferrari, Martin and Pimentel introduced the {\sl competition
interface} in the last-passage picture. This is a path
$k\mapsto \vp_k\in\Zb_+^2$ ($k\in\Zb_+$),
defined as a function of $\{G_{ij}\}$:  first $\vp_0=(0,0)$,
and then for $k\geq 0$
\be
\vp_{k+1}= \begin{cases}
\vp_k+(1,0) &\text{if $G(\vp_k+(1,0))<G(\vp_k+(0,1))$,}\\
\vp_k+(0,1) &\text{if $G(\vp_k+(1,0))>G(\vp_k+(0,1))$.}
\end{cases}\label{eq:def-vphi}
\ee
In other words, $\vp$ takes up-right steps, always choosing
the smaller of the two possible $G$-values.

The term ``competition interface'' is justified by the following picture.
Instead of having the unit squares centered at the integer points
as in Figure \ref{fig:iniri},  draw the squares so that
their corners coincide with integer points.
Label the  squares  by their northeast
corners, so that  the square
$(i-1,i]\times(j-1,j]$ is labeled  the
$(i,j)$-square. Regard  the last-passage time $G_{ij}$
as the time when the $(i,j)$-square becomes occupied.  Color the square $(0,0)$
white. Every other square gets either a red or a blue color:
 squares   to
the left and above the path $\vp$ are colored red, and squares   to
the right and below  $\vp$ blue.  Then the red squares are those whose
maximal path $\wh\pi$ passes through $(0,1)$, while the blue squares are
 those whose
maximal path $\wh\pi$ passes through $(1,0)$.  These can be regarded as
two competing ``infections'' on the $(i,j)$-plane, and $\vp$ is the
interface between them.

The competition interface represents the evolution of a second-class
particle, and macroscopically it follows the characteristics.
This was one of the main points for \cite{compint}.  In the
present setting the competition interface
 is the time reversal
of the maximal path $\wh\pi$, as we explain more precisely
in Section \ref{sec:prel}  below. This
connection allows us to establish the order of the transversal
fluctuations of the competition interface in the equilibrium setting.
To put this in precise notation, we introduce
\be
v(n)=\inf\{i\,:\,(i,n)=\vp_k\text{ for some }k\geq 0\}
\label{eq:def-vn}
\ee
and
\[
w(m)=\inf\{j\,:\,(m,j)=\vp_k\text{ for some }k\geq 0\}
\]
with the usual convention  $\inf\emptyset=\infty$. In other words,
$(v(n),n)$ is the leftmost point of the competition interface on
the horizontal line $j=n$, while $(m,w(m))$ is the lowest such
point on the vertical line $i=m$.  They are connected by the
implication
\be
v(n)\geq m\ \Longrightarrow \ w(m)<n
\label{eq:v-w-ineq}\ee
as can be seen from a picture.
Transposition $\om_{ij}\mapsto\om_{ji}$
of the $\om$-array interchanges $v$ and $w$.

Given $m$ and $n$,
let
\be
Z^{*\vr}=[m-v(n)]^+-[n-w(m)]^+\label{eq:zstar}
\ee
denote the signed distance from the point $(m,\,n)$
to the point where $\vp_k$ first hits either of the lines $j=n$ ($Z^{*\vr}>0$) or $i=m$ ($Z^{*\vr}<0$). Precisely one of the two terms
 contributes to the difference.
When we let $m=m(t)$ and $n=n(t)$ according to \eqref{eq:char-mn}, we have the $t$-dependent version $Z^{*\vr}(t)$.
Time reversal will show that in distribution
$Z^{*\vr}(t)$ is equal to $Z^\vr(t)$. (The notation $Z^*$ is used in anticipation of this time reversal connection.) Consequently
\begin{cor}
Theorem \ref{tm:Z-1} is true word for word when $Z^\vr(t)$ is replaced by $Z^{*\vr}(t)$.
\end{cor}
%
%
%


\subsection{Results for the rarefaction fan}\label{sc:rareres}

We now partially generalize the previous results to
arbitrary boundary conditions that are bounded by
the equilibrium boundary conditions of \eqref{eq:bondis}.
Let  $\{\om_{ij}\}$ be distributed  as in \eqref{eq:bondis}.
Let  $\{\hat{\om}_{ij}\}$ be another array
defined on the same probability space
 such that $\hat{\om}_{00}=0$,
 $\hat{\om}_{ij}=\om_{ij}$ for $i,j \geq 1$, and
\be
\hat{\om}_{i0}\leq \om_{i0} \quad\text{and}\quad \hat{\om}_{0j}\leq
\om_{0j} \quad \forall\ i,j\geq 1.\label{eq:rareom}
\ee
In particular, $\hat{\om}_{i0}=\hat{\om}_{0j}=0$ is admissible here.
Section \ref{sc:rarepart} below explains how
these  boundary conditions can represent the so-called
 rarefaction fan situation of simple exclusion.

Let  $\hat{G}(t)$ denote the weight of the maximal path to $(m,n)$ of
\eqref{eq:char-mn}, using the $\{\hat{\om}_{ij}\}$ array.
\begin{tm}\label{tm:rare}
Fix $0<\al <1$. There exists a
constant $C=C(\al, \vr)$ such that for all $t\ge1$ and $a>0$,
\[
\Pv \{ |\hat{G}(t)-t| > at^{1/3}\} \leq Ca^{-3\al/2}.
\]
\end{tm}
Define also $\hat{Z}_l(t)$ as the $i$-coordinate of the right-most point on the horizontal line $j=l$ of the right-most maximal path to $(m,n)$, and $\hat{Y}_l(t)$ as the $i$-coordinate of the left-most point on the horizontal line $j=l$ of the left-most maximal path to $(m,n)$. (In this general setting we no longer necessarily have a unique maximizing path because we have not ruled out a dependence of $\{\hat{\om}_{i0},\hat{\om}_{0j}\}$ on $\{\hat{\om}_{ij}\}_{i,j \geq 1}$.)
\begin{tm}\label{tm:raretrans}
For all $0<\al<1$ there exists $C=C(\al,\,\vr)$, such that for all $a>0$, $s\le t$ with $t\ge1$ and $(k,l)=(\fl{(1-\vr)^2s},\fl{\vr^2s})$,
\[
\Pv\{\hat{Z}_l(t)\geq k+at^{2/3}\}\leq Ca^{-3\al},\quad\text{and}\quad\Pv\{\hat{Y}_l(t)\le k-at^{2/3}\}\leq Ca^{-3\al}.
\]
\end{tm}

\section{Particle systems and queues}\label{sc:part}

The proofs in our paper will only use the last-passage description of the model. However, we would like to point out several other pictures one can attach to the last-passage
model.  An immediate one is the
 totally asymmetric simple exclusion process (TASEP).
  The boundary conditions \eqref{eq:bondis}  of the last-passage
model   correspond to  TASEP  in  equilibrium,
as seen by a ``typical'' particle right after its jump.
We also briefly discuss queues, and an augmentation of the last-passage
picture that describes a deposition model with column growth,
as in \cite{fluct}.

\subsection{The totally asymmetric simple exclusion process}
\label{sec:tasep}
This process
describes particles that jump unit steps  to the right
 on the integer lattice $\Zb$,  subject
to the exclusion rule that permits at most one particle per site.
The state of the process is a $\{0,1\}$-valued
sequence
$\wt{\un\eta}=\{\wt\eta_x\}_{x\in\Zb}$,
with the interpretation that  $\wt\eta_x=1$ means that site
$x$ is occupied by a particle, and
 $\wt\eta_x=0$  that $x$ is vacant.
The dynamics of the process are such that each $(1,0)$ pair
in the state
becomes a $(0,1)$ pair at rate 1, independently of the rest
of the state. In other words, each particle jumps to a vacant
site on its  right at rate 1, independently of other particles.
 The extreme points of the set of spatially
 translation-invariant equilibrium   distributions
 of this
process are the Bernoulli($\vr$) distributions $\nu^\vr$ indexed
by particle density  $0\leq\vr\leq1$.  Under $\nu^\vr$ the
occupation variables $\{\wt\eta_x\}$ are i.i.d.\ with
mean $\Ev^{\vr}(\wt\eta_x)=\vr$.

The Palm distribution of a particle system describes the equilibrium
distribution as seen from a ``typical'' particle. For a function $f$ of
$\wt{\un\eta}$, the Palm-expectation is
\[
\wh\Ev^\vr(f(\wt{\un\eta}))=\frac{\Ev^\vr(f(\wt{\un\eta})\cdot\wt\eta_0)}
{\Ev^\vr(\wt\eta_0)}
\]
in terms of the equilibrium expectation, see e.g.\ Port and Stone
\cite{posto}. Due to $\wt\eta_x\in\{0,\,1\}$,  for TASEP the Palm distribution is
the original Bernoulli($\vr$)-equilibrium conditioned on $\wt\eta_0=1$.

\begin{tm}[Burke]
Let $\wt{\un\eta}$ be a totally asymmetric simple exclusion process started from
the Palm distribution (i.e.\ a particle at the origin, Bernoulli measure elsewhere).
Then the position of the particle started at the origin is marginally a Poisson
process with jump rate $1-\vr$.
\end{tm}
The theorem follows from considering the inter-particle distances as M/M/1 queues.
Each of these distances is geometrically distributed, which is the stationary
distribution for the corresponding queue. Departure processes from these queues,
which correspond to TASEP particle jumps, are marginally Poisson due to Burke's
Theorem for queues, see e.g.\ Br\'emaud \cite{bremaud} for details. The Palm
distribution is important in this argument, as selecting a ``typical''
TASEP-particle assures that the inter-particle distances (or the lengths of the
queues) are geometrically distributed. For instance, the first particle to the
left of the origin in an ordinary Bernoulli equilibrium will \emph{not} see a
geometric distance to the next particle on its right.

Shortly we will explain how the boundary conditions
\eqref{eq:bondis} correspond to
 TASEP
started from Bernoulli($\vr$) measure, \emph{conditioned on $\wt\eta_0(0)=0$
and\linebreak [4]$\wt\eta_1(0)=1$}, i.e.\ a hole at the origin and a particle at site one
initially.  It will be convenient to give all particles and holes labels
that they retain as they jump (particles to the right, holes to the
left).
The particle initially at site one is labeled $P_0$, and the hole
initially at the origin is  labeled $H_0$. After this, all particles
are labeled with integers \emph{from right to left}, and all holes
\emph{from
left to right}. The position of particle $P_j$ at time $t$ is $P_j(t)$,
and the position of hole $H_i$  at time $t$ is $H_i(t)$.
Thus initially
\begin{align*}
  \dotsm < P_{3}(0)< P_{2}(0)&<P_{1}(0) < H_0(0)=0\\
& <1=P_0(0)
< H_1(0) < H_2(0)< H_3(0) <\dotsm
\end{align*}
  Since particles never
jump over each other, $P_{j+1}(t)<P_j(t)$ holds  at all times $t\geq 0$,
and by the same  token  also
$H_i(t)<H_{i+1}(t)$.

It turns out that this perturbation of the Palm distribution does
not entirely spoil Burke's Theorem.
\begin{cor}\label{cr:burke}
Marginally, $P_0(t)-1$ and $-H_0(t)$
are two independent Poisson processes with
respective jump rates $1-\vr$ and $\vr$.
\end{cor}
\begin{proof}
The evolution of $P_0(t)$ depends only on the initial
configuration\linebreak[4] $\{\wt\eta_x(0)\}_{x>1}$ and the Poisson
clocks governing the jumps over the edges $\{x\to x+1\}_{x\geq1}$. The
evolution of $H_0(t)$ depends only on the initial configuration $\{\wt\eta_x(0)\}_{x<0}$
and the Poisson clocks governing the jumps over the edges $\{x\to x+1\}_{x<0}$.
Hence $P_0(t)$ and $H_0(t)$ are independent. Moreover, $\{\wt\eta_x(0)\}_{x>1,\ x<0}$
is Bernoulli($\vr$) distributed, just like in the Palm distribution. Hence
Burke's Theorem applies to $P_0(t)$. As for $H_0(t)$, notice that $\un1-\un{\wt\eta}(t)$,
with $\un1_{x}\equiv1$, is a TASEP with holes and particles interchanged and
particles jumping to the left. Hence Burke's Theorem applies to $-H_0(t)$.
\end{proof}%

Now we can state the precise connection with the last-passage model.
For $i,j\geq 0$ let $T_{ij}$ denote the time when particle $P_j$ and hole $H_i$ exchange
places, with $T_{00}=0$.  Then
\[
\text{the processes $\{G_{ij}\}_{i,j\geq0}$ and $\{T_{ij}\}_{i,j\geq0}$ are equal in
distribution.  }
\]
For the marginal distributions on the
$i$- and $j$-axes  we see the truth
of the statement from Corollary \ref{cr:burke}.  More generally,
we can compare  the growing cluster
\[
\mathcal C(t)=\{(i,j)\in\Zb_+^2: T_{ij}\leq t\}
\]
with $\Ac(t)$   defined by \eqref{eq:def-Ac},
and observe that they are  countable state Markov chains
with the same  initial state and  identical bounded jump rates.

Since each particle jump corresponds to exchanging places
with a particular hole, one can deduce that at time $T_{ij}$,
\be
P_j(T_{ij})=i-j+1 \quad\text{ and }\quad H_i(T_{ij})=i-j.\label{eq:places}
\ee



By
the queuing interpretation of the
TASEP, we represent particles as ser\-vers, and the holes between $P_j$ and $P_{j-1}$
as customers in the queue of server $j$. Then \emph{the occupation of the
last-passage point $(i,\,j)$ is the same event as the
completion of the service of customer $i$
by server $j$}. This infinite system of queues is equivalent to a constant rate totally asymmetric zero range process.

\subsection{The rarefaction fan}\label{sc:rarepart}

The classical rarefaction fan initial condition for TASEP
is constructed with two densities $\la_\ell>\la_r$.
Initially  particles to the left of the origin
obey Bernoulli $\la_\ell$ distributions, and particles to the
right of the origin follow Bernoulli $\la_r$ distributions.
Of interest here is the behavior of a second-class particle
or the competition interface, and we refer the reader to
articles \cite{serf, compint, fermarpim, mountguiol, seck}

Following the development of the previous section,
 condition this initial
 measure on having a hole $H_0$ at $0$, and a particle $P_0$
at $1$. Then as  observed earlier, $H_0$
jumps to the left according to a Poisson$(\la_\ell)$ process, while
$P_0$ jumps to the right according to a Poisson$(1-\la_r)$ process.
To represent this situation in the last-passage picture,  choose
boundary weights $\{\hat{\om}_{i0}\}$ i.i.d.~Exp($1-\la_r$), and
 $\{\hat{\om}_{0j}\}$ i.i.d.~Exp($\la_\ell$),  corresponding to the waiting
times of $H_0$ and $P_0$. Suppose  $\la_\ell>\vr>\la_r$ and
$\om$ is the $\vr$-equilibrium
boundary condition defined by \eqref{eq:bondis}.
Then we have the
stochastic domination  $\om_{i0}\ge \hat{\om}_{i0}$ and
 $\om_{0j}\ge \hat\om_{0j}$, and we can realize these inequalities
by coupling the boundary weights.   The proofs
of Section \ref{sc:gen-bd} show that in fact one need not
insist on exponential boundary weights $\{\hat{\om}_{i0}, \hat{\om}_{0j}\}$,
but instead only inequality \eqref{eq:rareom} is required for
the fluctuations.

\subsection{A deposition model}

In this section we describe a deposition model that gives a direct graphical connection between
the
TASEP and the last-passage percolation.  This point of view
is not needed for the later proofs, hence we only give a
brief explanation.

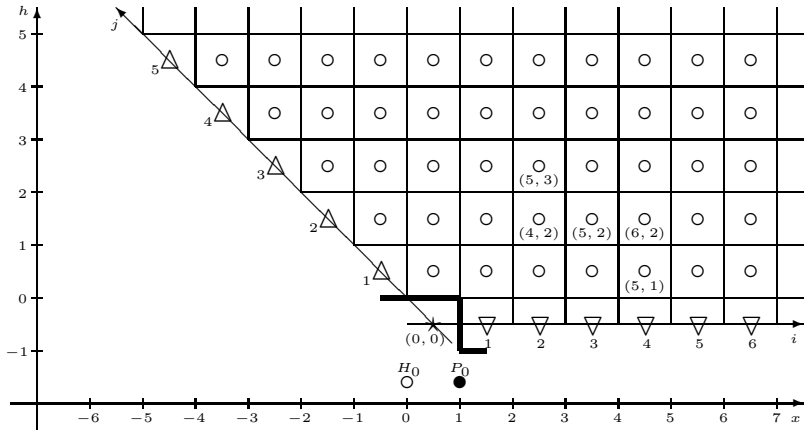
\begin{figure}[ht]
\begin{center}
\begin{picture}(320,180)(-20,-20)
\put(-10,0){\vector(1,0){300}}
\put(0,-10){\vector(0,1){160}}
\put(140,30){\vector(1,0){150}}
\put(157,23){\vector(-1,1){127}}
\put(140,40){\line(1,0){150}}
\put(120,60){\line(1,0){170}}
\put(100,80){\line(1,0){190}}
\put(80,100){\line(1,0){210}}
\put(60,120){\line(1,0){230}}
\put(40,140){\line(1,0){250}}
\put(40,140){\line(0,1){10}}
\put(60,120){\line(0,1){30}}
\put(80,100){\line(0,1){50}}
\put(100,80){\line(0,1){70}}
\put(120,60){\line(0,1){90}}
\put(140,40){\line(0,1){110}}
\multiput(160,30)(20,0){7}{\line(0,1){120}}
\put(147,27){\large$\star$}
\multiput(166,26)(20,0){6}{\large$\triangledown$}
\multiput(126,47)(-20,20){5}{\large$\vartriangle$}
\multiput(147,47)(20,0){7}{\large$\circ$}
\multiput(127,67)(20,0){8}{\large$\circ$}
\multiput(107,87)(20,0){9}{\large$\circ$}
\multiput(87,107)(20,0){10}{\large$\circ$}
\multiput(67,127)(20,0){11}{\large$\circ$}
\put(139,23){\tiny$(0,0)$}
\put(285,23){\tiny$i$}
\put(28,142){\tiny$j$}
\put(169,21){\tiny$1$}
\put(189,21){\tiny$2$}
\put(209,21){\tiny$3$}
\put(229,21){\tiny$4$}
\put(249,21){\tiny$5$}
\put(269,21){\tiny$6$}
\put(123,45){\tiny$1$}
\put(103,65){\tiny$2$}
\put(83,85){\tiny$3$}
\put(63,105){\tiny$4$}
\put(43,125){\tiny$5$}
\put(202,63){\tiny$(5,2)$}
\put(222,43){\tiny$(5,1)$}
\put(182,83){\tiny$(5,3)$}
\put(182,63){\tiny$(4,2)$}
\put(222,63){\tiny$(6,2)$}
\multiput(20,-2)(20,0){14}{\line(0,1){4}}
\multiput(-2,20)(0,20){7}{\line(1,0){4}}
\put(285,-7){\tiny$x$}
\put(-7,147){\tiny$h$}
\put(15,-7){\tiny$-6$}
\put(35,-7){\tiny$-5$}
\put(55,-7){\tiny$-4$}
\put(75,-7){\tiny$-3$}
\put(95,-7){\tiny$-2$}
\put(115,-7){\tiny$-1$}
\put(138,-7){\tiny$0$}
\put(158,-7){\tiny$1$}
\put(178,-7){\tiny$2$}
\put(198,-7){\tiny$3$}
\put(218,-7){\tiny$4$}
\put(238,-7){\tiny$5$}
\put(258,-7){\tiny$6$}
\put(278,-7){\tiny$7$}
\put(-12,18){\tiny$-1$}
\put(-7,38){\tiny$0$}
\put(-7,58){\tiny$1$}
\put(-7,78){\tiny$2$}
\put(-7,98){\tiny$3$}
\put(-7,118){\tiny$4$}
\put(-7,138){\tiny$5$}
\put(137,5){\large$\circ$}
\put(157,5){\large$\bullet$}
\put(136,12){\tiny$H_0$}
\put(156,12){\tiny$P_0$}
{\linethickness{2pt}
\put(130,40){\line(1,0){30}}
\put(160,20){\line(0,1){20}}
\put(160,20){\line(1,0){10}}
}
\end{picture}
\end{center}
\caption{The initial configuration}\label{fig:ini}
\end{figure}

We start by tilting the $j$-axis and all the vertical columns of Figure \ref{fig:iniri} by 45 degrees, resulting in Figure \ref{fig:ini}.
This picture represents
the same initial situation as Figure \ref{fig:iniri}, but note
that now the $j$-coordinates must be read  in the direction $\nwarrow$.
(As before, some  squares are labeled with their $(i,j)$-coordinates.)
The $i-j$ tilted coordinate system is embedded in an $x-h$ orthogonal system.

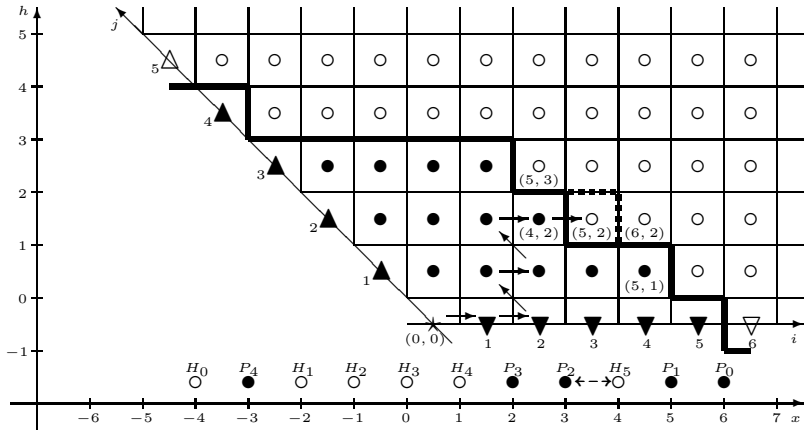
\begin{figure}[ht]
\begin{center}
\begin{picture}(320,180)(-20,-20)
\put(-10,0){\vector(1,0){300}}
\put(0,-10){\vector(0,1){160}}
\put(140,30){\vector(1,0){150}}
\put(157,23){\vector(-1,1){127}}
\put(140,40){\line(1,0){150}}
\put(120,60){\line(1,0){170}}
\put(100,80){\line(1,0){190}}
\put(80,100){\line(1,0){210}}
\put(60,120){\line(1,0){230}}
\put(40,140){\line(1,0){250}}
\put(40,140){\line(0,1){10}}
\put(60,120){\line(0,1){30}}
\put(80,100){\line(0,1){50}}
\put(100,80){\line(0,1){70}}
\put(120,60){\line(0,1){90}}
\put(140,40){\line(0,1){110}}
\multiput(160,30)(20,0){7}{\line(0,1){120}}
\put(147,27){\large$\star$}
\multiput(166,26)(20,0){5}{\large$\blacktriangledown$}
\put(266,26){\large$\triangledown$}
\multiput(126,47)(-20,20){4}{\large$\blacktriangle$}
\put(46,127){\large$\vartriangle$}
\multiput(247,47)(20,0){2}{\large$\circ$}
\multiput(207,67)(20,0){4}{\large$\circ$}
\multiput(187,87)(20,0){5}{\large$\circ$}
\multiput(87,107)(20,0){10}{\large$\circ$}
\multiput(67,127)(20,0){11}{\large$\circ$}
\multiput(127,67)(-20,20){2}{\multiput(0,0)(20,0){4}{\large$\bullet$}}
\multiput(147,47)(20,0){5}{\large$\bullet$}
\put(139,23){\tiny$(0,0)$}
\put(285,23){\tiny$i$}
\put(28,142){\tiny$j$}
\put(169,21){\tiny$1$}
\put(189,21){\tiny$2$}
\put(209,21){\tiny$3$}
\put(229,21){\tiny$4$}
\put(249,21){\tiny$5$}
\put(269,21){\tiny$6$}
\put(123,45){\tiny$1$}
\put(103,65){\tiny$2$}
\put(83,85){\tiny$3$}
\put(63,105){\tiny$4$}
\put(43,125){\tiny$5$}
\multiput(20,-2)(20,0){14}{\line(0,1){4}}
\multiput(-2,20)(0,20){7}{\line(1,0){4}}
\put(285,-7){\tiny$x$}
\put(-7,147){\tiny$h$}
\put(15,-7){\tiny$-6$}
\put(35,-7){\tiny$-5$}
\put(55,-7){\tiny$-4$}
\put(75,-7){\tiny$-3$}
\put(95,-7){\tiny$-2$}
\put(115,-7){\tiny$-1$}
\put(138,-7){\tiny$0$}
\put(158,-7){\tiny$1$}
\put(178,-7){\tiny$2$}
\put(198,-7){\tiny$3$}
\put(218,-7){\tiny$4$}
\put(238,-7){\tiny$5$}
\put(258,-7){\tiny$6$}
\put(278,-7){\tiny$7$}
\put(-12,18){\tiny$-1$}
\put(-7,38){\tiny$0$}
\put(-7,58){\tiny$1$}
\put(-7,78){\tiny$2$}
\put(-7,98){\tiny$3$}
\put(-7,118){\tiny$4$}
\put(-7,138){\tiny$5$}
\put(202,63){\tiny$(5,2)$}
\put(222,43){\tiny$(5,1)$}
\put(182,83){\tiny$(5,3)$}
\put(182,63){\tiny$(4,2)$}
\put(222,63){\tiny$(6,2)$}
\multiput(57,5)(160,0){2}{\large$\circ$}
\multiput(97,5)(20,0){4}{\large$\circ$}
\put(77,5){\large$\bullet$}
\multiput(177,5)(60,0){2}{\multiput(0,0)(20,0){2}{\large$\bullet$}}
\put(203,6){$\dashrightarrow$}
\put(203,6){$\dashleftarrow$}
\put(56,12){\tiny$H_0$}
\put(96,12){\tiny$H_1$}
\put(116,12){\tiny$H_2$}
\put(136,12){\tiny$H_3$}
\put(156,12){\tiny$H_4$}
\put(216,12){\tiny$H_5$}
\put(76,12){\tiny$P_4$}
\put(176,12){\tiny$P_3$}
\put(196,12){\tiny$P_2$}
\put(236,12){\tiny$P_1$}
\put(256,12){\tiny$P_0$}
\multiput(155,33)(20,37){2}{
\multiput(0,0)(20,0){2}{\vector(1,0){11}}}
\multiput(185,35)(0,20){2}{\vector(-1,1){10}}
\put(175,50){\vector(1,0){11}}
{\linethickness{2pt}
\put(50,120){\line(1,0){30}}
\put(80,100){\line(1,0){100}}
\multiput(180,80)(60,-40){2}{\line(1,0){20}}
\put(200,60){\line(1,0){40}}
\put(260,20){\line(1,0){10}}
\multiput(80,100)(160,-60){2}{\line(0,1){20}}
\multiput(180,80)(20,-20){2}{\line(0,1){20}}
\put(260,20){\line(0,1){20}}
\multiput(201,80)(4,0){5}{\line(1,0){2}}
\multiput(220,61)(0,4){5}{\line(0,1){2}}
}
\end{picture}
\end{center}
\caption{A possible move at a later time}\label{fig:later}
\end{figure}

Figure \ref{fig:later} shows the later situation that corresponds
to Figure \ref{fig:laterri}. As before, the thickset line is the boundary of the squares belonging to $\Ac(t)$ of \eqref{eq:def-Ac}. Whenever it makes sense, the height $h_x$ of a column $x$ is defined as the $h$-coordinate (i.e.\ the vertical height) of the thickset line above the edge $[x,\,x+1]$ on the $x$-axis. Define the increments $\eta_x=h_{x-1}-h_x$ and notice that, whenever defined, $\eta_x\in\{0,\,1\}$ due to the tilting we made. The last passage rules, converted for this picture, tell us that occupation of a new square happens at rate one unless it would violate $\eta_x\in\{0,\,1\}$ for some $x$. Moreover, one can read that the occupation of a square $(i,\,j)$ is the same event as the pair $(\eta_{i-j},\,\eta_{i-j+1})$ changing from $(1,\,0)$ to $(0,\,1)$. Comparing this to \eqref{eq:places} leads us to the conclusion that $\eta_x$, whenever defined, is the occupation variable of the simple exclusion process that corresponds to the last passage model. This way one can also conveniently include the particles ($\eta_x=1$) and holes ($\eta_x=0$) on the $x$-axis, as seen on the figures. Notice also that the time-increment $h_x(t)-h_x(0)$ is the cumulative particle current across the bond $[x,\,x+1]$.

\subsection{The characteristics}\label{sec:char}

One-dimensional conservative particle systems have the conservation law
\[
\pt_t\vr(t,\,x)+\pt_xf(\vr(t,\,x))=0
\]
under the Eulerian hydrodynamic scaling, where $\vr(t,\,x)$ is the expected particle number per site and $f(\vr(t,\,x))$ is the macroscopic particle flux around the rescaled position $x$ at the rescaled time $t$, see e.g.\ \cite{cl} for details. Disturbances of the solution propagate with the characteristic speed $f'(\vr)$.
The macroscopic particle flux for TASEP is $f(\vr)=\vr(1-\vr)$,
and consequently the characteristic speed is $f'(\vr)=1-2\vr$.
Thus the characteristic curve started
at the origin  is  $t\mapsto (1-2\vr)t$.  To identify the  point $(m,n)$
in the last-passage picture that corresponds to th%
is
curve, we reason  approximately. Namely, we look for $m$ and $n$ such that hole $H_m$ and particle $P_n$ interchange positions at around time $t$ and the characteristic position $(1-2\vr)t$. By time $t$, that particle $P_n$ has jumped over approximately $(1-\vr)t$ sites due to Burke's Theorem. Hence at time zero, $P_n$ is approximately at position $(1-2\vr)t-(1-\vr)t=-\vr t$. Since the particle density is $\vr$, the particle labels around this position are $n\approx\vr^2t$ at time zero. Similarly, holes travel at a speed $-\vr$, so hole $H_m$ starts from approximately $(1-2\vr)t+\vr t$. They have density $1-\vr$, which indicates $m\approx(1-\vr)^2t$.
Thus we are led to consider the point
$(m,n)=(\tfl{(1-\vr)^2t}, \tfl{\vr^2t})$ as done in \eqref{eq:char-mn}.

\section{Preliminaries}
\label{sec:prel}

We turn to establish some basic facts and tools.
First an extension
of  Corollary \ref{cr:burke}
to show that Burke's Theorem holds for every hole and particle
 in the last-passage picture. Define
\[
\ba
I_{ij}:&=G_{ij}-G_{\{i-1\}j}\quad\text{for }i\geq1,\ j\geq0,\quad\text{ and}\\
J_{ij}:&=G_{ij}-G_{i\{j-1\}}\quad\text{for }i\geq0,\ j\geq1.
\ea
\]
$I_{ij}$ is the time it takes for particle $P_j$ to jump again after its jump to site $i-j$.
$J_{ij}$ is the time it takes for hole $H_i$ to jump again after its jump to site $i-j+1$.
Applying the last passage rules \eqref{eq:grec} shows
\be
\begin{split}
I_{ij}&=G_{ij}-G_{\{i-1\}j}\\
&=(G_{\{i-1\}j}\lor G_{i\{j-1\}})+\om_{ij}-G_{\{i-1\}\{j-1\}}-(G_{\{i-1\}j}-G_{\{i-1\}\{j-1\}})\label{eq:iij}\\
&=(J_{\{i-1\}j}\lor I_{i\{j-1\}})+\om_{ij}-J_{\{i-1\}j}\\
&=(I_{i\{j-1\}}-J_{\{i-1\}j})^++\om_{ij}.
\end{split}
\ee
Similarly,
\be
J_{ij}=(J_{\{i-1\}j}-I_{i\{j-1\}})^++\om_{ij}.\label{eq:jij}
\ee
For later use, we define
\be
X_{\{i-1\}\{j-1\}}=I_{i\{j-1\}}\land J_{\{i-1\}j}.\label{eq:xdef}
\ee
\begin{lm}\label{lm:markov}
Fix $i,\,j\geq1$. If $I_{i\{j-1\}}$ and $J_{\{i-1\}j}$ are independent
exponentials with respective parameters $1-\vr$ and $\vr$, then $I_{ij}$, $J_{ij}$, and $X_{\{i-1\}\{j-1\}}$ are jointly independent exponentials with respective parameters $1-\vr$, $\vr$, and $1$.
\end{lm}
\begin{proof}
As the variables $I_{i\{j-1\}}$, $J_{\{i-1\}j}$ and $\om_{ij}$ are independent, we use \eqref{eq:iij}, \eqref{eq:jij} and \eqref{eq:xdef} to write the joint moment generating function as
\begin{multline*}
M_{I_{ij},\,J_{ij},\,X_{\{i-1\}\{j-1\}}}(s,\,t,\,u):=\Ev\e{sI_{ij}+tJ_{ij}+uX_{\{i-1\}\{j-1\}}}\\
=\Ev\e{s(I_{i\{j-1\}}-J_{\{i-1\}j})^++t(J_{\{i-1\}j}-I_{i\{j-1\}})^++u(I_{i\{j-1\}}\land J_{\{i-1\}j})}\cdot\Ev\e{(s+t)\om_{ij}}
\end{multline*}
where it is defined. Then, with the assumption of the lemma
and the definition of $\om_{ij}$, elementary calculations show
\[
M_{I_{ij},\,J_{ij},\,X_{\{i-1\}\{j-1\}}}(s,\,t,\,u)=\frac{\vr\cdot(1-\vr)}{(1-\vr-s)\cdot(\vr-t)\cdot(1-u)}.
\]
\end{proof}

Let $\Sigma$ be the set of doubly-infinite down-right paths
in the first quadrant of the $(i,j)$-coordinate system.
 In terms of the sequence of points visited
 a path $\si\in\Sigma$ is given by
\[
\si=\{\dots\to
 (p_{-1}, \,q_{-1})\to (p_0,\,q_0)\to(p_1,\,q_1)\to\dots
\to(p_l,\,q_l)\to\dots\}
\]
with all $p_l,q_l\geq 0$ and
 steps
\[
(p_{l+1},\,q_{l+1})-(p_l,\,q_l)= 
\begin{cases} (1,0) &\text{(direction $\to$ in Figure \ref{fig:iniri}), or}\\
(0,-1) &\text{(direction $\downarrow$  in Figure \ref{fig:iniri}).}
\end{cases}
\]
The interior
of the set enclosed by $\si$ is defined by
\[
\mathcal B(\si)=\{(i,\,j)\,:\,0\leq i<p_l,\ 0\leq j<q_l\text{ for some }(p_l,\,q_l)\in\si\}.
\]
The last-passage time increments along $\si$ are the variables
\[
 Z_l(\si)=G_{p_{l+1}q_{l+1}}-G_{p_lq_l}=
 \left\{\ba
&I_{p_{l+1}q_{l+1}},&&\text{if }(p_{l+1},\,q_{l+1})-(p_l,\,q_l)=(1,0),\\
&J_{p_lq_l},&&\text{if }(p_{l+1},\,q_{l+1})-(p_l,\,q_l)=(0,-1),
\ea\right.
\]
for $l\in\Zb$.
We admit the possibility that  $\si$ is the union of the $i$- and
$j$-coordinate
axes, in which case $\mathcal B(\si)$ is empty.
\begin{lm}\label{lm:NE}
For any $\si\in\Sigma$,  the random variables
\be
\bigl\{\{X_{ij}\,:\,(i,\,j)\in\mathcal B(\si)\},\  \{Z_l(\si): l\in\Zb\}
\bigr\}\label{eq:mind}
\ee
are mutually independent, $I$'s with $\Expd(1-\vr)$, $J$'s with
 $\Expd(\vr)$, and $X$'s with $\Expd(1)$ distribution.
\end{lm}
\begin{proof}  We first consider the countable set of
paths that join the $j$-axis to the $i$-axis, in other words
those for which there exist
 finite $n_0<n_1$ such that
$p_{n}=0$ for $n\leq n_0$ and  $q_{n}=0$ for $n\geq n_1$.
For these paths we argue by induction on $\mathcal B(\si)$.
When $\mathcal B(\si)$ is the empty set,
 the statement reduces to the independence of $\om$-values on the $i$- and $j$-axes
which is part of the set-up.

Now given an arbitrary $\si\in\Sigma$ that connects the $j$- and the $i$-axes, consider
a {\sl growth corner} $(i,j)$ for $\mathcal B(\si)$,
by which we mean that for some index $l\in\Zb$,
\[ (p_{l-1}, q_{l-1}), (p_{l}, q_{l}), (p_{l+1}, q_{l+1})
\;=\;
(i,j+1), (i,j), (i+1,j). \]
A new valid $\wt\si\in\Sigma$ can be produced by
replacing the above points with
\[ (\wt p_{l-1}, \wt q_{l-1}), (\wt p_{l}, \wt q_{l}),
(\wt p_{l+1}, \wt q_{l+1})
\;=\;
(i,j+1), (i+1,j+1), (i+1,j) \]
and now $\mathcal B(\wt\si)=\mathcal B(\si)\cup\{(i,j)\}$.

The change inflicted on the set of random variables
\eqref{eq:mind} is that
\be
\{I_{\{i+1\}j},\ J_{i\{j+1\}}\}\label{eq:rem}
\ee
has been replaced by
\be
\{I_{\{i+1\}\{j+1\}},\ J_{\{i+1\}\{j+1\}},\ X_{ij}\}.\label{eq:add}
\ee
By \eqref{eq:iij}--\eqref{eq:jij}
variables  \eqref{eq:add} are determined by \eqref{eq:rem} and
$\om_{\{i+1\}\{j+1\}}$.  If we assume
inductively that  $\si$ satisfies the conclusion we seek, then
so does $\wt\si$ by
 Lemma \ref{lm:markov} and because in the situation under
consideration  $\om_{\{i+1\}\{j+1\}}$ is
independent of the variables in \eqref{eq:mind}.

 For an arbitrary $\sigma$ the statement follows
because the independence of the random variables in \eqref{eq:mind}
follows from independence of finite subcollections.
Consider any square $R=\{0\leq i,j\leq M\}$  large enough so that
the corner $(M,M)$ lies outside $\si\cup\mathcal B(\si)$.
Then the $X$- and $Z(\sigma)$-variables associated to $\si$
that lie in $R$ are a subset of the variables of a certain
path $\wt\si$ that goes through the points $(0,M)$ and $(M,0)$.
Thus the variables
in \eqref{eq:mind} that lie inside an arbitrarily large square
are independent.
\end{proof}


By applying Lemma \ref{lm:NE}
to a path that contains the horizontal line  $q_l\equiv j$ we get
a version of Burke's theorem: particle
$P_j$ obeys a Poisson process after time
$G_{0j}$ when it  ``enters the last-passage picture.''
The vertical line  $p_l\equiv i$
gives the corresponding  statement for hole $H_i$.

Example 2.10.2 of Walrand \cite{walrand} gives an
intuitive understanding of this result. Our initial state
corresponds to the situation when particle $P_0$
and hole $H_0$ have just exchanged places
 in an equilibrium system of queues. $H_0$ is
therefore a customer who has just moved from queue 0 to queue 1.
By that Example, this customer sees an equilibrium system of
queues every time he jumps. Similarly, any new customer arriving
to the queue of particle $P_1$ sees an equilibrium queue system in
front, so Burke's theorem extends to the region between $P_0$ and $H_0$.

Up-right turns do not have independence:
variables $I_{ij}$ and $J_{i\{j+1\}}$, or $J_{ij}$ and $I_{\{i+1\}j}$
are not independent.

The same inductive argument with a growing cluster $\mathcal B(\si)$
 proves a result that corresponds to a coupling of
 two exclusion systems $\eta$ and $\wt\eta$ where the
latter has a higher density of particles.
However, the lemma is a purely deterministic
statement.

\begin{lm}\label{lm:basic1}
Consider two assignments of values $\{\om_{ij}\}$ and $\{\wt\om_{ij}\}$
that satisfy $\om_{00}=\wt\om_{00}=0$,
 $\om_{0j}\geq\wt\om_{0j}$, $\om_{i0}\leq\wt\om_{i0}$,
and $\om_{ij}=\wt\om_{ij}$  for all $i,j\geq 1$.     Then all
increments satisfy  $I_{ij}\leq\wt I_{ij}$ and  $J_{ij}\geq\wt J_{ij}$.
\end{lm}
\begin{proof}  One  proves by induction that the statement holds for
all increments between points in $\si\cup\mathcal B(\si)$ for
those paths $\si\in\Sigma$ for which $\mathcal B(\si)$ is finite.
If $\mathcal B(\si)$ is empty the statement is the assumption
made on the $\om$- and $\wt\om$-values on the $i$- and $j$-axes. The induction
step that adds a growth corner to $\mathcal B(\si)$
 follows from equations \eqref{eq:iij} and \eqref{eq:jij}.
\end{proof}

\subsection{The reversed process.}

Fix $m>0$ and $n>0$, and define
\[
H_{ij}=G_{mn}-G_{ij}
\]
for $0\leq i\leq m$, $0\leq j\leq n$. This is the time needed to ``free'' the
point $(i,\,j)$ in the reversed process, started from the moment when $(m,\,n)$
becomes occupied.  For $0\leq i<m$ and $0\leq j<n$,
\[
\ba
H_{ij}&
=-((G_{\{i+1\}j}-G_{mn})\land(G_{i\{j+1\}}-G_{mn}))\\
&\qquad+((G_{\{i+1\}j}-G_{ij})\land(G_{i\{j+1\}}-G_{ij}))\\
&=H_{\{i+1\}j}\lor H_{i\{j+1\}}+X_{ij}
\ea
\]
with definition \eqref{eq:xdef} of the $X$-variables. Taking this and Lemma \ref{lm:NE} into account, we see that the $H$-process is a
copy of the original $G$-process, but with reversed coordinate directions.
Precisely speaking,  define $\om^*_{00}=0$,
and then for $0<i\leq m$, $0<j\leq n$: $\om^*_{i0}=I_{\{m-i+1\}n}$, $\om^*_{0j}=J_{m\{n-j+1\}}$,
and $\om^*_{ij}=X_{\{m-i\}\{n-j\}}$.  Then  $\{\om^*_{ij}: 0\leq i\leq m, 0\leq j\leq n\}$
is distributed like $\{\om_{ij}: 0\leq i\leq m, 0\leq j\leq n\}$ in \eqref{eq:bondis}, and
 the process
\be
G^*_{ij}=H_{\{m-i\}\{n-j\}}\label{eq:def-K}
\ee
for $0\leq i\leq m$, $0\leq j\leq n$ satisfies
\[
G^*_{ij}=(G^*_{\{i-1\}j}\lor G^*_{i\{j-1\}})+\om^*_{ij}, \qquad 0\leq i\leq m\,,\, 0\leq j\leq n
\]
(with the formal assumption $G^*_{\{-1\}j}=G^*_{i\{-1\}}=0$), see \eqref{eq:grec}.
Thus the pair  $(G^*,\om^*)$ has the same distribution as $(G,\om)$ in a fixed
rectangle $\{0\leq i\leq m\}\times \{ 0\leq j\leq n\}$. Throughout the paper quantities defined in the reversed process will be denoted by a superscript $^*$, and they will always be equal in distribution to their original forward versions.

\subsection{Exit point and competition interface.}

For integers $x$ define
\be
U^\vr_x=G_{\{x^+\}\{x^-\}}=\begin{cases}
\sum_{i=0}^x \om_{i0}, &x\geq 0\\
\sum_{j=0}^{-x} \om_{0j}, &x\leq 0.
\end{cases} \label{eq:udef}
\ee
 Referring to the two
coordinate systems in Figure \ref{fig:ini},
this is  the last-passage time of
 the point on the $(i,j)$-axes  above point  $x$ on the $x$-axis.
This point is on the $i$-axis
if $x\geq0$ and  on the $j$-axis if $x\leq0$.

Fix  integers $m\geq x^+\lor 1$, $n\geq x^-\lor 1$, and define $\Pi_x(m,n)$
as the set of directed paths $\pi$ connecting
$(x^+\lor 1,\,x^-\lor 1)$ and $(m,\,n)$
using allowable steps \eqref{eq:allow}. Then let
\be
A_x=A_x(m,n)=\max_{\pi\in\Pi_x(m,n)}\sum_{(p,q)\in\pi}\om_{pq}\label{eq:adef}
\ee
be the maximal weight collected by a path from
$(x^+,\,x^-)$ to $(m,\,n)$ that immediately exits the axes,
and does not  count $\om_{x^+x^-}$.
Notice that $A_{-1}=A_0=A_1$ and this value is the last-passage
time from $(1,1)$ to $(m,n)$ that completely ignores the
boundaries, or in other words, sets the boundary values
$\om_{i0}$ and $\om_{0j}$  equal
to zero.

By the continuity of the exponential distribution
 there is an a.s.\ unique path $\wh\pi$
 from $(0,\,0)$ to $(m,\,n)$ which
collects the maximal weight $G^\vr_{mn}$.   Earlier we defined
the {\sl exit point}  $Z^\vr\in\Zb$ to represent  the last
point of this path on either the $i$-axis or the $j$-axis.
Equivalently we can now state that
 $Z^\vr$ is the a.s.\ unique integer for which
\[
G_{mn}^\vr=U^\vr_{Z^\vr}+A_{Z^\vr}.
\]
Simply because the maximal path $\wh\pi$ necessarily goes through either
$(0,1)$ or $(1,0)$, $Z^\vr$ is always nonzero.


Recall the definition of the competition interface in \eqref{eq:def-vphi}.
Now we can observe
that the competition interface is the time reversal
of the maximal path $\wh\pi$. Namely, the competition interface
of the reversed process follows the maximal path $\wh\pi$ backwards
from the corner  $(m,n)$, until it hits either the $i$- or the
$j$-axis.  To make a precise statement, let us represent the
a.s.\ unique maximal last-passage
path, with exit point $Z^\vr$ as defined above, as
\[
\wh\pi=\{ (0,0)=\wh\pi_0 \to \wh\pi_1 \to \dotsm \to \wh\pi_{\lvert Z^\vr\rvert}
\to\dots \to \wh\pi_{m+n}=(m,n)\},
\]
where  $\{\wh\pi_{\lvert Z^\vr\rvert+1}
\to\dots \to \wh\pi_{m+n}\}$ is the portion of the path
that  resides in the interior
$\{1,\dotsc,m\}\times\{1,\dotsc,n\}$.

\begin{lm}
  Let $\vp^*$
be the competition interface constructed for the process
$G^*$ defined by \rm{\eqref{eq:def-K}}. Then
$\vp_k^*=(m,n)-\wh\pi_{m+n-k}$ for $0\leq k\leq m+n-\lvert Z^\vr\rvert$.
\end{lm}
\begin{proof}  Starting from $\wh\pi_{m+n}=(m,n)$,
the maximal path $\wh\pi$ can be constructed backwards
step by step  by always moving to the maximizing point
of the right-hand side of \eqref{eq:grec}.
This is the same as constructing the competition interface
for the reversed process $G^*$ by \eqref{eq:def-vphi}. Since $G^*$ is not constructed
outside the rectangle $\{0,\dotsc,m\}\times\{0,\dotsc,n\}$, we
cannot assert what the competition interface does after the
point
\[\vp^*_{m+n-\lvert Z^\vr\rvert}=
\wh\pi_{\lvert Z^\vr\rvert}=\begin{cases} (Z^\vr, 0) &\text{if $Z^\vr>0$}\\
(0,-Z^\vr) &\text{if $Z^\vr<0$.}
\end{cases}
\qedhere
\]
\end{proof}
Notice that, due to this lemma, $Z^{*\vr}$ defined in \eqref{eq:zstar} is indeed $Z^\vr$ defined in the reversed process, which justifies the argument following \eqref{eq:zstar}.

The competition interface
bounds  the regions  where the boundary conditions on the  axes are
felt.  From this we can get useful
 bounds between last-passage times under different
boundary conditions.  This is  the last-passage model equivalent
of the common use of second-class particles to control discrepancies
between coupled  interacting particle systems. In the next lemma,
the superscript $\Wc$ represents the west boundary ($j$-axis) of the
$(i,j)$-plane. Remember that $(v(n),n)$ is the left-most point of the competition interface on the horizontal line $j=n$
computed in terms of the $G$-process (see \eqref{eq:def-vn}).

\begin{lm} Let $G^{\Wc=0}$ be the last-passage times of a system where
we set $\om_{0j}=0$ for all $j\geq 1$. Then for $v(n)<m_1<m_2$,
\begin{align*}
A_0(m_2,n)-A_0(m_1,n)&\leq G^{\Wc=0}(m_2,n)-G^{\Wc=0}(m_1,n)\\
&= G(m_2,n)-G(m_1,n).
\end{align*}
\label{lm:coup-3}
\end{lm}
\begin{proof} The first inequality is a consequence of
Lemma \ref{lm:basic1}, because computing $A_0$ is the same
as computing $G$ with all boundary values $\om_{i0}=\om_{0j}=0$ (and in fact this inequality is valid for all $m_1<m_2$.)
The equality $G(m,n)=G^{\Wc=0}(m,n)$ for $m>v(n)$ follows because
the maximal path $\wh\pi$ for $G(m,n)$ goes through $(1,0)$ and
hence does not see the boundary values $\om_{0j}$. Thus this
same path $\wh\pi$ is maximal for $G^{\Wc=0}(m,n)$ too.
\end{proof}

If we set $\om_{i0}=0$ (the south boundary denoted by $\Sc$)
instead, we get this
statement: for $0\leq m_1<m_2\leq v(n)$,
\be
\begin{split}
A_0(m_2,n)-A_0(m_1,n)&\geq G^{\Sc=0}(m_2,n)-G^{\Sc=0}(m_1,n)\\
&= G(m_2,n)-G(m_1,n).
\label{eq:coup-4}
\end{split}
\ee

\subsection{A coupling on the $i$-axis}

Let $1>\la>\vr>0$. As a common realization of the exponential weights $\om_{i0}^\la$ of $\Expd(1-\la)$ and $\om_{i0}^\vr$ of $\Expd(1-\vr)$ distribution, we write
\be
\om_{i0}^\la=\frac{1-\vr}{1-\la}\cdot\om_{i0}^\vr.\label{eq:cou}
\ee
We will use this coupling later for different purposes. We will also need
\[
\Vv(\om_{i0}^\la-\om_{i0}^\vr)=\left(\frac{1-\vr}{1-\la}-1\right)^2\cdot\frac{1}{(1-\vr)^2}=\left(\frac{1}{1-\la}-\frac{1}{1-\vr}\right)^2.
\]

\subsection{Exit point and the variance of the last-passage time}

With these preliminaries
we can prove  the key lemma that links the variance of the last-passage
time to the weight collected along the axes.

\begin{lm}\label{lm:ns}
Fix $m,\,n$ positive integers. Then
\be
\ba
\Vv(G^\vr_{mn})&=\frac{n}{\vr^2}-\frac{m}{(1-\vr)^2}+\frac{2}{1-\vr}\cdot\Ev(U_{Z^{\vr+}}^\vr)\\
&=\frac{m}{(1-\vr)^2}-\frac{n}{\vr^2}+\frac{2}{\vr}\cdot\Ev(U_{-Z^{\vr-}}^\vr),
\ea\label{eq:ns}
\ee
where $Z^\vr$ is the a.s.\ unique exit point of the maximal path from $(0,\,0)$ to $(m,\,n)$.
\end{lm}
\begin{proof}
We label the total increments along the sides of the rectangle by compass
directions:
\[
\Wc=G_{0n}^\vr-G_{00}^\vr,\quad\Nc=G_{mn}^\vr-G_{0n}^\vr,\quad\Ec=G_{mn}^\vr-G_{m0}^\vr,\quad\Sc=G_{m0}^\vr-G_{00}^\vr.
\]
As $\Nc$ and $\Ec$ are independent by Lemma \ref{lm:NE}, we have
\be
\ba
\Vv(G^\vr_{mn})&=\Vv(\Wc+\Nc)\\
&=\Vv(\Wc)+\Vv(\Nc)+2\Cov(\Sc+\Ec-\Nc,\,\Nc)\\
&=\Vv(\Wc)-\Vv(\Nc)+2\Cov(\Sc,\,\Nc).
\ea\label{eq:varg}
\ee
We now modify the $\om$-values in $\Sc$. Let $\la=\vr+\ve$ and
apply \eqref{eq:cou}, without changing the other values
 $\{\om_{ij}: i\geq 0,j\geq1\}$.  Quantities of the altered last-passage
 model will be marked with a superscript $\ve$. In
this new process, $\Sc^\ve$ has a Gamma$(m,\,1-\vr-\ve)$ distribution with density
\[
f_\ve(s)=\frac{(1-\vr-\ve)^m\cdot\e{-(1-\vr-\ve)s}\cdot s^{m-1}}{(m-1)!}
\]
for $s>0$, whose $\ve$-derivative is
\be
\pt_\ve f_\ve(s)=sf_\ve(s)-\frac{m}{1-\vr-\ve}\cdot f_\ve(s).\label{eq:fe}
\ee
Given the sum $\Sc^\ve$, the joint distribution of $\{\om_{i0}\}_{1\leq i\leq m}$ is independent of the parameter $\ve$, hence the quantity $\Ev(\Nc^\ve\,|\,\Sc^\ve=s)=\Ev(\Nc\,|\,\Sc=s)$ does not depend on $\ve$. Therefore, using \eqref{eq:fe} we have
\be
\begin{split}
\pt_\ve\Ev(\Nc^\ve)\Bigr|_{\ve=0}&=\pt_\ve\int_0^\infty\Ev(\Nc\,|\,\Sc=s)f_\ve(s)\,\di s\Bigr|_{\ve=0}\\
&=\int_0^\infty\Ev(\Nc\,|\,\Sc=s)\cdot s\cdot f_0(s)\,\di s-\frac{m}{1-\vr}\int_0^\infty\Ev(\Nc\,|\,\Sc=s)f_0(s)\,\di s\\
&=\Ev(\Nc\Sc)-\frac{m}{1-\vr}\cdot\Ev(\Nc)=\Cov(\Nc,\,\Sc).\label{eq:dcov}
\end{split}
\ee
Next we compute the same quantity by a different approach. Let $Z$ and $Z^\ve$ be the exit points of the maximal paths to $(m,\,n)$ in the original and the modified processes, respectively. Similarly, $U_x$ and $U^\ve_x$ are the weights as defined by \eqref{eq:udef} for the two processes. Hence $U_Z$ is the weight collected on the $i$ or $j$ axis by the maximal path of the original process. Then
\[
\ba
\Nc^\ve-\Nc&=(\Nc^\ve-\Nc)\cdot{\bf1}\{Z^\ve=Z\}+(\Nc^\ve-\Nc)\cdot{\bf1}\{Z^\ve\neq Z\}\\
&=(U_Z^\ve-U_Z)\cdot{\bf1}\{Z^\ve=Z\}+(\Nc^\ve-\Nc)\cdot{\bf1}\{Z^\ve\neq Z\}\\
&=(U_Z^\ve-U_Z)+(\Nc^\ve-\Nc-U_Z^\ve+U_Z)\cdot{\bf1}\{Z^\ve\neq Z\}.
\ea
\]
As $\om$ values are only changed on the $i$-axis, the first term is rewritten as
\[
U_Z^\ve-U_Z=U_{Z^+}^\ve-U_{Z^+}=\left(\frac{1-\vr}{1-\vr-\ve}-1\right)U_{Z^+}=\frac{\ve}{1-\vr-\ve}\cdot U_{Z^+}
\]
by \eqref{eq:cou}. We show that the expectation
of the second term is $\mathfrak o(\ve)$.   Note that the
increase $\Nc^\ve-\Nc$ is bounded by $\Sc^\ve-\Sc$. Hence
\be
\begin{split}
&\Ev[(\Nc^\ve-\Nc-U_Z^\ve+U_Z)\cdot{\bf1}\{Z^\ve\neq Z\}]\\
&\qquad \leq\Ev[(\Nc^\ve-\Nc)\cdot{\bf1}\{Z^\ve\neq Z\}]
\;\leq\;\Ev[(\Sc^\ve-\Sc)\cdot{\bf1}\{Z^\ve\neq Z\}]\\
&\qquad \leq\bigl(\Ev[(\Sc^\ve-\Sc)^2]\bigr)^\frac12\cdot
\bigl(\Pv\{Z^\ve\neq Z\}\bigr)^\frac12.
\end{split}\label{eq:err}
\ee
To show that the probability is of the order of $\ve$, notice that the exit point of the maximal path can only differ in the modified process from the one of the original process, if for some $Z<k\leq m$, $U^\ve_k+A_k>U^\ve_Z+A_Z$ with $Z$ of the original process (see \eqref{eq:adef} for the definition of $A_i$). Therefore,
\[
\ba
\Pv\{Z^\ve\neq Z\}&=\Pv\{U^\ve_k+A_k>U^\ve_Z+A_Z\text{ for some }Z<k\leq m\}\\
&=\Pv\{U^\ve_k-U^\ve_Z>A_Z-A_k\text{ for some }Z<k\leq m\}\\
&=\Pv\{U^\ve_k-U^\ve_Z>A_Z-A_k\geq U_k-U_Z\text{ for some }Z<k\leq m\}\\
&\leq\Pv\{U^\ve_k-U^\ve_i>A_i-A_k\geq U_k-U_i\text{ for some }0\leq i<k\leq m\}\\
&\leq\sum_{0\leq i<k\leq m}\Pv\{U^\ve_k-U^\ve_i>A_i-A_k\geq U_k-U_i\}.
\ea
\]
We also used the definition of $Z$ in the third equality, via $A_Z+U_Z\geq A_k+U_k$. Notice that $A$'s and $U$'s are independent for fixed indices. Hence with $\mu$ denoting the distribution of $A_i-A_k$, we write
\begin{multline*}
\Pv\{U^\ve_k-U^\ve_i>A_i-A_k\geq U_k-U_i\}\\
\ba
&=\int\Pv\{U^\ve_k-U^\ve_i>x\geq U_k-U_i\}\di\mu(x)\\
&\leq\sup_x\Pv\{U^\ve_k-U^\ve_i>x\geq U_k-U_i\}\\
&=\sup_x\Pv\Bigl\{\frac{1-\vr}{1-\vr-\ve}\cdot(U_k-U_i)>x\geq U_k-U_i\Bigr\}\\
&=\sup_x\Pv\Bigl\{x\geq U_k-U_i>x\Bigl(1-\frac{\ve}{1-\vr}\Bigr)\Bigr\}.
\ea
\end{multline*}
Since $U_k-U_i$ has a Gamma distribution, the supremum above is $\mathcal O(\ve)$, which shows the bound on $\Pv\{Z^\ve\neq Z\}$. The first factor on the right-hand side of \eqref{eq:err},
\[
\bigl(\Ev[(\Sc^\ve-\Sc)^2]\bigr)^{1/2}
=\frac{\ve}{1-\vr-\ve}\cdot\bigl(\Ev[\Sc^2]\bigr)^{1/2},
\]
is of order $\ve$.
Hence the error term \eqref{eq:err} is $\mathfrak o(\ve)$, and we conclude
\[
\pt_\ve\Ev(\Nc^\ve)\Bigr|_{\ve=0}=\frac{1}{1-\vr}\cdot\Ev(U_{Z^+}).
\]
The proof of the first statement is then completed by this display, \eqref{eq:varg} and \eqref{eq:dcov}, as $\Wc$ and $\Nc$ are Gamma-distributed by Lemma \ref{lm:NE}. The second statement follows in a similar way, using $\Cov(\Wc,\,\Ec)$.
\end{proof}

\begin{lm}\label{lm:vvlr}
Let $0<\vr\leq\la<1$. Then
\[
\Vv(G_{mn}^\la)\leq\frac{\vr^2}{\la^2}\cdot\Vv(G^\vr_{mn})+m\cdot\left(\frac{1}{(1-\la)^2}-\frac{\vr^2}{\la^2(1-\vr)^2}\right).
\]
\end{lm}
\begin{proof}
The proof is based on the coupling described by \eqref{eq:cou}, and a similar one $\om_{0j}^\la=\frac\vr\la\cdot\om_{0j}^\vr$ on the $j$ axis. Note that in this coupling, when changing from $\vr$ to $\la$, we are increasing the weights on the $i$-axis and decreasing the weights on the $j$-axis, which clearly implies $Z^\vr\leq Z^\la$. Also, we remain in the stationary situation, so \eqref{eq:ns} remains valid for $\la$. As $U_{-x^-}^\vr$ is non-increasing in $x$, this implies
\[
U^\la_{-Z^{\la-}}=\frac\vr\la\cdot U^\vr_{-Z^{\la-}}\leq\frac\vr\la\cdot U^\vr_{-Z^{\vr-}}.
\]
We substitute this into the second line of \eqref{eq:ns} to get
\[
\ba
\Vv(G^\la_{mn})&=\frac{m}{(1-\la)^2}-\frac{n}{\la^2}+\frac{2}{\la}\cdot\Ev(U_{-Z^{\la-}}^\la)\\
&\leq\frac{m}{(1-\la)^2}-\frac{n}{\la^2}+\frac{2\vr}{\la^2}\cdot\Ev(U_{-Z^{\vr-}}^\vr)\\
&=\frac{\vr^2}{\la^2}\cdot\left(\frac{m}{(1-\vr)^2}-\frac{n}{\vr^2}+\frac{2}{\vr}\cdot\Ev(U_{-Z^{\vr-}}^\vr)\right)\\
&\qquad+m\cdot\left(\frac{1}{(1-\la)^2}-\frac{\vr^2}{\la^2(1-\vr)^2}\right)\\
&=\frac{\vr^2}{\la^2}\cdot\Vv(G^\vr_{mn})+m\cdot\left(\frac{1}{(1-\la)^2}-\frac{\vr^2}{\la^2(1-\vr)^2}\right).
\ea
\]
\end{proof}


\section{Upper bound}
\label{sc:ub}

We turn  to proving the upper bounds in Theorems \ref{tm:G-var-1} and \ref{tm:Z-1}.
We have a fixed  density
$\vr\in(0,1)$, and to study the last-passage times $G^\vr$
 along the  characteristic,
we define the
dimensions of the last-passage rectangle as
\be
m(t)=\fl{(1-\vr)^2t}\qquad\text{and}\qquad n(t)=\fl{\vr^2t}\label{eq:mndef}
\ee
with a parameter $t\to\infty$.
The quantities $A_x$, $Z$ and $G_{mn}$ connected to these indices
are denoted by $A_x(t)$, $Z(t)$, $G(t)$.
In the proofs we  need to consider different  boundary conditions
 \eqref{eq:bondis}  with $\vr$ replaced by $\la$.
This will be indicated by a superscript.
 However,  the
superscript $\la$ only changes the boundary conditions
and not the dimensions  $m(t)$ and $n(t)$, always defined by
\eqref{eq:mndef} with a fixed $\vr$.
 Moreover,
we apply the coupling \eqref{eq:cou} on the $i$-axis and $\om_{0j}^\la=\frac\vr\la\cdot\om_{0j}^\vr$ on the $j$-axis.
The weights $\{\om_{ij}\}_{i,\,j\geq1}$ in the interior
will not be affected by changes in boundary conditions,
so in particular  $A_x(t)$ will not either. Since
$G^\la(t)$ chooses the  maximal path,
\[
U^\la_z+A_z(t)\leq G^\la(t)
\]
for all $1\leq z\leq m(t)$
and all densities $0<\la<1$.
Consequently, for integers $u\geq0$ and densities $\la\geq\vr$,
\be
\ba
\Pv\{Z^\vr(t)>u\}&=\Pv\{\exists z>u\,:\,U^\vr_z+A_z(t)=G^\vr(t)\}\\
&\leq\Pv\{\exists z>u\,:\,U^\vr_z-U^\la_z+G^\la(t)\geq G^\vr(t)\}\\
&=\Pv\{\exists z>u\,:\,U^\la_z-U^\vr_z\leq G^\la(t)-G^\vr(t)\}\\
&\leq\Pv\{U^\la_u-U^\vr_u\leq G^\la(t)-G^\vr(t)\}.
\ea\label{eq:pzu}
\ee
The last step is justified by $\la\geq\vr$ and the coupling
\eqref{eq:cou}.
Set
\be
\la_u=\frac{\vr}{\sqrt{(1-\vr)^2-{u}/{t}}+\vr}.\label{eq:lu}
\ee
This density maximizes
\[
\Ev(U^\la_u)-\Ev(G^\la(t))=\frac{u}{1-\la}
-\frac{\fl{(1-\vr)^2t}}{1-\la}-\frac{\fl{\vr^2t}}{\la}
\]
if the integer parts are dropped. The expectation
$\Ev(G^\la_{mn})$ is
computed as  $\Ev(G^\la_{0n})+\Ev(G^\la_{mn}-G^\la_{0n})$
with the help of Lemma \ref{lm:NE}.
 Some useful identities for future computations:
\be
\la_u\geq\vr,\quad\frac{1}{\la_u}=1+\frac{\sqrt{(1-\vr)^2-{u}/{t}}}{\vr},\quad\frac{1}{1-\la_u}=1+\frac{\vr}{\sqrt{(1-\vr)^2-{u}/{t}}}.\label{eq:usef}
\ee

\begin{lm}
With  $0\leq u\leq (1-\vr)^2t$ and  $\la_u$ of \eqref{eq:lu},
\begin{multline*}
\Ev(U^{\la_u}_u-U^\vr_u-G^{\la_u}(t)+G^\vr(t))\\
\geq
\frac{t\vr}{1-\vr}\left((1-\vr)-\sqrt{(1-\vr)^2-{u}/{t}}\,\right)^2
-\frac{u/t}{\vr(1-\vr)}.
\end{multline*}
\end{lm}
\begin{proof}
By Lemma \ref{lm:NE} and \eqref{eq:mndef}
\begin{align*}
&\Ev(U^{\la_u}_u-U^\vr_u-G^{\la_u}(t)+G^\vr(t))\\
&\quad =\frac{u}{1-\la_u}-\frac{u}{1-\vr}-\frac{\fl{(1-\vr)^2t}}{1-\la_u}+
\frac{\fl{(1-\vr)^2t}}{1-\vr}-\frac{\fl{\vr^2t}}{\la_u}+
\frac{\fl{\vr^2t}}{\vr}.
\end{align*}
First we remove the integer parts.
Since $\la_u\geq\vr$,
\[
-\frac{\fl{(1-\vr)^2t}}{1-\la_u}+
\frac{\fl{(1-\vr)^2t}}{1-\vr} \geq
-\frac{(1-\vr)^2t}{1-\la_u}+
\frac{(1-\vr)^2t}{1-\vr}.
\]
For the other integer parts
\begin{align*}
-\frac{\fl{\vr^2t}}{\la_u}+
\frac{\fl{\vr^2t}}{\vr} &\geq  -\frac{\vr^2t}{\la_u}+
\frac{\vr^2t}{\vr}  -\frac1\vr+\frac1\la_u\\
&=-\frac{\vr^2t}{\la_u}+
\frac{\vr^2t}{\vr}-\,\frac{1-\vr-\sqrt{(1-\vr)^2-u/t}}{\vr}\\
&\geq -\frac{\vr^2t}{\la_u}+
\frac{\vr^2t}{\vr} -\frac{u/t}{\vr(1-\vr)}
\end{align*}
The last term above
is  the last term of the bound in the statement of the lemma.
It remains to check that after  the integer parts have been
removed from the mean, the remaining quantity  equals the main
 term of the bound.
\begin{multline*}
\frac{u}{1-\la_u}-\frac{u}{1-\vr}-\frac{(1-\vr)^2t}{1-\la_u}+\frac{(1-\vr)^2t}{1-\vr}-\frac{\vr^2t}{\la_u}+\frac{\vr^2t}{\vr}\\
\ba
&=[u-(1-\vr)^2t]\cdot\left[1+\frac{\vr}{\sqrt{(1-\vr)^2-{u}/{t}}}-\frac{1}{1-\vr}\right]\\
&\qquad-\vr^2t\cdot\left[1+\frac{\sqrt{(1-\vr)^2-{u}/{t}}}{\vr}-\frac{1}{\vr}\right]\\
&=t\cdot\frac{\vr}{1-\vr}\left((1-\vr)-\sqrt{(1-\vr)^2-{u}/{t}}\right)^2.
\ea
\end{multline*}
\end{proof}

\begin{lm}\label{lm:middle}
For any $8\vr^{-2}(1-\vr)^2\leq u\leq(1-\vr)^2t$,
\[
\Ev(U^{\la_u}_u-U^\vr_u-G^{\la_u}(t)+G^\vr(t))\geq
\frac{\vr}{8(1-\vr)^3}\cdot\frac{u^2}{t}.
\]
\end{lm}
\begin{proof}
Assumption $u\geq 8\vr^{-2}(1-\vr)^2$ implies that the last term
of the bound from the previous lemma satisfies
\[
-\,\frac{u/t}{\vr(1-\vr)}\geq -\,\frac{\vr}{8(1-\vr)^3}\cdot\frac{u^2}{t}.
\]
Thus it remains to prove
\[
\left((1-\vr)-\sqrt{(1-\vr)^2-{u}/{t}}\,\right)^2\geq\frac{1}{4(1-\vr)^2}\cdot\frac{u^2}{t^2}.
\]
This is easy to check in the form
\[
\left(C-\sqrt{C^2-x}\right)^2
\geq\frac{1}{4C^2}\cdot x^2
,
\]
where  $x=u/t$, $C=1-\vr$ and then $x\leq C^2$.
\end{proof}

\begin{lm}\label{lm:vg}
For any $0\leq u\leq\frac34(1-\vr)^2t$,
\[
\Vv(G^{\la_u}(t)-G^\vr(t))\leq\frac{8}{1-\vr}\cdot\Ev(U^\vr_{Z^\vr(t)^+})+
\frac{8(u+1)}{(1-\vr)^2}
\]
\end{lm}
\begin{proof}
We start with substituting \eqref{eq:mndef} into Lemma \ref{lm:vvlr}
(integer parts can  be dropped without violating the inequality):
\[
\Vv(G^{\la_u}(t))\leq\frac{\vr^2}{\la_u^2}\cdot\Vv(G^\vr(t))+t\cdot\left(\frac{(1-\vr)^2}{(1-\la_u)^2}-\frac{\vr^2}{\la_u^2}\right).
\]
Utilizing  \eqref{eq:usef},
\[
\frac{(1-\vr)^2}{(1-\la_u)^2}-\frac{\vr^2}{\la_u^2}
=\left(\sqrt{(1-\vr)^2-{u}/{t}}+\vr\right)^2\cdot\frac{{u}/{t}}{(1-\vr)^2-{u}/{t}}.
\]
Since the expression in parentheses is not larger than 1,
${u}/{t}\leq\frac34(1-\vr)^2$, and $\vr\leq\la_u$, it follows that
\[
\Vv(G^{\la_u}(t))\leq\Vv(G^\vr(t))+\frac{4}{(1-\vr)^2}\cdot u.
\]
Then we proceed with Lemma \ref{lm:ns} and \eqref{eq:mndef}:
\begin{align*}
&\Vv(G^{\la_u}(t)-G^\vr(t))\leq2\Vv(G^{\la_u}(t))+2\Vv(G^{\vr}(t))\\
&\qquad \leq4\Vv(G^{\vr}(t))+\frac{8}{(1-\vr)^2}\cdot u\\
&\qquad =\frac{8}{1-\vr}\cdot\Ev(U^\vr_{Z^\vr(t)^+})+4\frac{\fl{\vr^2t}}{\vr^2}
-4\frac{\fl{(1-\vr)^2t}}{(1-\vr)^2}
+\frac{8}{(1-\vr)^2}\cdot u\\
&\qquad\leq\frac{8}{1-\vr}\cdot\Ev(U^\vr_{Z^\vr(t)^+})+\frac{8(u+1)}{(1-\vr)^2}
 \qedhere
\end{align*}
\end{proof}

\begin{lm}\label{lm:vu}
With the application of the coupling \eqref{eq:cou}, for any $0\leq u\leq\frac34(1-\vr)^2t$  
we have
\[
\Vv(U_u^{\la_u}-U_u^\vr)\leq u\cdot\frac{\vr^2}{(1-\vr)^2}.
\]
\end{lm}
\begin{proof}
By that coupling,
\[
\Vv[U_u^{\la_u}-U_u^\vr]=\Vv\left[\left(\frac{1-\vr}{1-\la_u}-1\right)U_u^\vr\right]=u\cdot\left(\frac{1-\vr}{1-\la_u}-1\right)^2\cdot\frac{1}{(1-\vr)^2},
\]
as $U^\vr_u$ is the sum of $u$ many independent $\Expd$$(1-\vr)$ weights.
Write
\[
\ba
\left(\frac{1-\vr}{1-\la_u}-1\right)\cdot\frac{1}{(1-\vr)}&=\frac{\sqrt{(1-\vr)^2-{u}/{t}}+\vr}{\sqrt{(1-\vr)^2-{u}/{t}}}-\frac{1}{(1-\vr)}\\
&\leq\frac{\frac12(1+\vr)}{\frac12(1-\vr)}-\frac{1}{(1-\vr)}=\frac{\vr}{1-\vr}.
\ea
\]
\end{proof}

\noindent
After these preparations, we continue the main argument from \eqref{eq:pzu}.
\begin{lm}\label{lm:zbou}
There exists a constant $C_1=C_1(\vr)$ such that
for any $u\geq\linebreak[4]8\vr^{-2}(1-\vr)^2$ and $t>0$,
\[
\Pv\{Z^\vr(t)>u\}\leq C_1\Bigl( \frac{t^2}{u^4}\cdot\Ev(U^\vr_{Z^\vr(t)^+})
+\frac{t^2}{u^3}\Bigr).
\]
\end{lm}
\begin{proof}
If $8\vr^{-2}(1-\vr)^2\leq u\leq(1-\vr)^2t$, then continuing
from \eqref{eq:pzu} and taking Lemma \ref{lm:middle} into account, we write
\[
\ba
\Pv\{Z^\vr(t)>u\}
&\leq\Pv\biggl\{U^{\la_u}_u-U^\vr_u\leq\Ev(U^{\la_u}_u-U^\vr_u)-
\frac{\vr}{16(1-\vr)^3}\cdot\frac{u^2}{t}\biggr\}\\
&\ +\Pv\biggl\{G^{\la_u}(t)-G^\vr(t)\geq\Ev(G^{\la_u}(t)-G^\vr(t))+
\frac{\vr}{16(1-\vr)^3}\cdot\frac{u^2}{t}\biggr\}\\
&\leq\Vv(U^{\la_u}_u-U^\vr_u)\cdot\frac{16^2(1-\vr)^6}{\vr^2}\cdot\frac{t^2}{u^4}\\
&\ +\Vv(G^{\la_u}(t)-G^\vr(t))\cdot\frac{16^2(1-\vr)^6}{\vr^2}\cdot\frac{t^2}{u^4}
\ea
\]
by Chebyshev's inequality.  If
$8\vr^{-2}(1-\vr)^2\leq u\leq\frac34(1-\vr)^2t$,
use Lemmas \ref{lm:vu} and \ref{lm:vg} to conclude
\begin{align*}
\Pv\{Z^\vr(t)>u\}&\leq 16^2(1-\vr)^4\cdot\frac{t^2}{u^3}
+8\cdot 16^2\cdot\frac{(1-\vr)^5}{\vr^2}\cdot\frac{t^2}{u^4}
\cdot\Ev(U^\vr_{Z^\vr(t)^+})\\
&+8\cdot 16^2\cdot\frac{(1-\vr)^4}{\vr^2}\cdot\frac{t^2(u+1)}{u^4}.
\end{align*}
When $\frac34(1-\vr)^2t<u\leq(1-\vr)^2t$, the previous display
 works for $\frac34u$. Hence by
\[
\Pv\{Z^\vr(t)>u\}\leq\Pv\{Z^\vr(t)> 3u/4\},
\]
the statement still holds, modified by a  factor of a  power of $4/3$.

Finally, the probability is trivially zero if $u>(1-\vr)^2t$.
\end{proof}

\hop
Fix a number $0<\al<1$, and define
\be
y=\frac{u}{\al(1-\vr)}.\label{eq:uy}
\ee

\begin{lm}\label{lm:ld}
We have the following large deviations estimate:
\[
\Pv\{U^\vr_u>y\}\leq\e{-(1-\vr)(1-\sqrt\al)^2y}.
\]
\end{lm}
\begin{proof}
We use the fact that $U^\vr_u=\sum\limits_{i=1}^u\om_{i0}$, where the $\om$'s are iid.\ $\Expd$$(1-\vr)$ variables. Fix $s$ with $1-\vr>s>0$. By the Markov inequality, we write
\[
\ba
\Pv\{U^\vr_u>y\}&=\Pv\{\e{sU^\vr_u}>\e{sy}\}\leq\e{-sy}\Ev(\e{sU^\vr_u})
=\e{-sy}\cdot\Bigl(\frac{1-\vr}{1-\vr-s}\Bigr)^u\\
&\leq\exp\Bigl(-sy+u\cdot\frac{s}{1-\vr-s}\Bigr).
\ea
\]
Substituting $u=\al(1-\vr)y$, the choice $s=(1-\vr)(1-\sqrt\al)$ minimizes the exponent, and yields the result.
\end{proof}

\begin{lm}
There exist finite
positive constants  $C_2=C_2(\al, \vr)$ and  $C_3=C_3(\al, \vr)$
 such that, for all
\[
r\geq\frac{8(1-\vr)}{\al\vr^2 \Ev(U^\vr_{Z^\vr(t)^+})}\,,
\]
we have the bound
\begin{multline*}
\Pv\{U^\vr_{Z^\vr(t)^+}>r\Ev(U^\vr_{Z^\vr(t)^+})\}\\
\leq \frac{C_2t^2}{[\Ev(U^\vr_{Z^\vr(t)^+})]^3}\cdot
\left(\frac{1}{r^3}+\frac{1}{r^4}\right)
+\exp\{-C_3 r\Ev(U^\vr_{Z^\vr(t)^+})\}.
\end{multline*}
\end{lm}
\begin{proof}
By  \eqref{eq:uy} and  Lemmas
\ref{lm:zbou} and \ref{lm:ld},  for
any $ y\geq 8\al^{-1}\vr^{-2}(1-\vr)$, and with an appropriately
defined new constant,
\begin{align*}
\Pv\{U^\vr_{Z^\vr(t)^+}>y\}&\leq\Pv\{Z^\vr(t)^+>u\}+\Pv\{U^\vr_u>y\}\\
&\leq
C_2 \Bigl(\,
\frac{t^2}{y^4}\cdot\Ev(U^\vr_{Z^\vr(t)^+})+\frac{t^2}{y^3}\Bigr)
+\e{-(1-\vr)(1-\sqrt\al)^2y}.
\end{align*}
Choose $y=r\Ev(U^\vr_{Z^\vr(t)^+})$.
\end{proof}

\begin{tm}
\[
\limsup_{t\to\infty}\frac{\Ev(U^\vr_{Z^\vr(t)^+})}{t^{2/3}}<\infty,
\quad\text{and}\quad\limsup_{t\to\infty}\frac{\Vv(G^\vr(t))}{t^{2/3}}<\infty.
\]
\label{tm:EUbd}\end{tm}
\begin{proof}
The first inequality implies the second one by Lemma \ref{lm:ns} and \eqref{eq:mndef}. To prove the first one, suppose that there exists a sequence $t_k\nearrow\infty$ such that
\[
\lim_{k\to\infty}\frac{\Ev(U^\vr_{Z^\vr(t_k)^+})}{t_k^{2/3}}=\infty.
\]
Then $\Ev(U^\vr_{Z^\vr(t_k)^+})>t_k^{2/3}$ for all large $k$'s,
and consequently
by the above lemma
\[
\Pv\{U^\vr_{Z^\vr(t_k)^+}>r\Ev(U^\vr_{Z^\vr(t_k)^+})\}\\
\leq C_2 \left(\frac{1}{r^3}+\frac{1}{r^4}\right)\frac{t_k^{2}}{[\Ev(U^\vr_{Z^\vr(t_k)^+})]^3}
+\exp(-C_3 r t_k^{2/3})
\]
for all $r\geq C_4 t_k^{-2/3}$.
This shows by dominated convergence that
\[
\int_0^\infty\Pv\{U^\vr_{Z^\vr(t_k)^+}>r
\Ev(U^\vr_{Z^\vr(t_k)^+})\}\di r\underset{k\to\infty}{\longrightarrow}0,
\]
which leads to the contradiction
\[
1 = \Ev\left(\frac{U^\vr_{Z^\vr(t_k)^+}}{\Ev(U^\vr_{Z^\vr(t_k)^+})}\right)\underset{k\to\infty}{\longrightarrow}0.
\]
\end{proof}

Combining Lemma \ref{lm:zbou} and Theorem \ref{tm:EUbd} gives
a tail bound on $Z$:

\begin{cor} Given any $t_0>0$
there exists a finite constant  $C_4=C_4(t_0, \vr)$
such that, for all $a>0$ and $t\geq t_0$,
\[
\Pv\{ Z^\vr(t)\geq at^{2/3}\} \leq C_4 a^{-3}.
\]
\label{cor:Z-ub}
\end{cor}

\section{Lower bound}
\label{sc:lb}
We abbreviate
 $(m,n)=\bigl( \fl{(1-\vr)^2t},\fl{\vr^2t}\bigr)$
throughout this section.

\begin{lm} Let $a,b>0$ be arbitrary positive numbers.
There exist finite  constants $t_0=t_0(a,b,\vr)$
and $C=C(\vr)$  such that, for all $t\geq t_0$,
\[
\Pv\Bigl\{ \,\sup_{1\leq z\leq \va t^{2/3}}\bigl(  U^\vr_z + A_z(t) -A_1(t)\bigr)
\geq bt^{1/3} \Bigr\} \leq C\va^3 (b^{-3}+b^{-6}).
\]
\label{lm:lb-1}
\end{lm}
\begin{proof}
The process
$\{U^\vr_z\}$ depends on the boundary $\{\om_{i0}\}$.
Pick a version\linebreak[4] $\{\om_{ij}\}_{1\leq i\leq m,1\leq j\leq n}$ of the  interior variables
independent of $\{\om_{i0}\}$.
If we use the reversed system
 \be\{\wt\om_{ij}=\om_{m-i+1,n-j+1}\}_{1\leq i\leq m,1\leq j\leq n}
\label{eq:reversal-1}
\ee
 to compute $A_z(m,n)$,
then this coincides with $A_1(m-z+1,n)$ computed with $\{\om_{ij}\}$.
Thus with this coupling (and some abuse of notation) we can
replace $A_z(m,n) -A_1(m,n)$ with $A_1(m-z+1,n)-A_1(m,n)$.
[Note that $A_1(m,n)$ is the same for $\om$ and $\wt\om$.]
Next pick a further independent version of boundary conditions
\eqref{eq:bondis}   with density $\la$. Use these and $\{\om_{ij}\}_{i,j\geq 1}$
to compute the last-passage times $G^\la$, together with a
 competition interface  $\vp^\la$ defined by
\eqref{eq:def-vphi} and the projections $v^\la$
defined by \eqref{eq:def-vn}. Then by  \eqref{eq:coup-4}, on the
event $v^\la(n)\geq m$,
 \[
A_1(m,n)-A_1(m-z+1,n)  \geq G^\la(m,n)-G^\la(m-z+1,n).
\]
 Set
\[
V^\la_z= G^\la(m,n)
-G^\la(m-z,n),
\]
a sum of $z$ i.i.d.\ Exp($1-\la$) variables. $V^\la$ is independent
of $U^\vr$. Combining these steps we get the bound
\be
\begin{split}
&\Pv\Bigl\{ \,\sup_{1\leq z\leq \va t^{2/3}}\bigl(  U^\vr_z + A_z(t) -A_1(t)\bigr)
\geq bt^{1/3} \Bigr\} \\
&\quad \leq
\Pv\bigl\{ v^\la(\fl{\vr^2t})<\fl{(1-\vr)^2t}\bigr\}
 + \Pv\Bigl\{ \,\sup_{1\leq z\leq \va t^{2/3}}(  U^\vr_z -V^\la_{z-1})
\geq bt^{1/3} \Bigr\}
\end{split}
\label{eq:lb-temp-2}
\ee

  Introduce a parameter $r>0$ whose value
will be specified later, and define \be
\la=\vr-rt^{-1/3}.\label{eq:def-temp-la} \ee
For the second probability on the right-hand
side of \eqref{eq:lb-temp-2}, define the martingale
$M_z=U^\vr_z -V^\la_{z-1}-\Ev(U^\vr_z -V^\la_{z-1})$, and note that
for $z\leq \va t^{2/3}$,
\begin{align*}
&\Ev(U^\vr_z -V^\la_{z-1})= \frac{z}{1-\vr}-\frac{z-1}{1-\la}
= \frac{zrt^{-1/3}}{(1-\vr)(1-\la)}+\frac1{1-\la} \\
&\leq  \frac{r\va t^{1/3}}{(1-\vr)^2} +\frac1{1-\vr}.
\end{align*}
As long as
\be
 b> r\va(1-\vr)^{-2}+t^{-1/3}(1-\vr)^{-1},
\label{eq:lb-temp-2-a}
\ee
we get by Doob's inequality, for any $p\geq 1$,
\be
\begin{split}
&\Pv\Bigl\{ \,\sup_{1\leq z\leq \va t^{2/3}}(  U^\vr_z -V^\la_{z-1})
\geq bt^{1/3} \Bigr\}\\
&\leq \Pv\Bigl\{ \,\sup_{1\leq z\leq \va t^{2/3}} M_z \geq  t^{1/3}\Bigl( b
 - \frac{r\va}{(1-\vr)^2} - \frac{t^{-1/3}}{1-\vr}\,\Bigr)\,\Bigr\}\\
&\leq \frac{C(p)  t^{-p/3}}{\bigl( b-r\va(1-\vr)^{-2}-t^{-1/3}(1-\vr)^{-1}\bigr)^p}
\Ev\bigl[\,\lvert  M_{\fl{\va t^{2/3}}} \rvert^p\,\bigr]\\
&\leq \frac{C(p,\vr)  \va^{p/2}}{\bigl( b-r\va(1-\vr)^{-2}-t^{-1/3}(1-\vr)^{-1}\bigr)^p}.
\end{split}
\label{eq:lb-temp-2-b}\ee
Now choose $t_0=4^3b^{-3}(1-\vr)^{-3}$. Then for $t\geq t_0$ the above
bound is dominated by
\[
\frac{C(p,\vr)  \va^{p/2}}
{\Bigl( \frac{3b}4 - \frac{r\va}{(1-\vr)^2}\Bigr)^p}
\]
which becomes $C(p,\vr)\va^{3} b^{-6}  $
once we choose
\be r=\frac{b(1-\vr)^2}{4\va}
\label{eq:def-temp-r-1}\ee
and $p=6$,
and change the constant $C(p,\vr)$.

For  the first probability on the right-hand
side of \eqref{eq:lb-temp-2},
introduce the time $s=(\vr/\la)^2t$. Then
\[
\Pv\bigl\{ v^\la(\fl{\vr^2t})<\fl{(1-\vr)^2t}\bigr\}=\Pv\bigl\{ v^\la(\fl{\la^2s})<\fl{\la^2(1-\vr)^2\vr^{-2}s}\bigr\}.
\]
Notice that since $\la<\vr$ here, $\fl{\la^2(1-\vr)^2\vr^{-2}s}\leq\fl{(1-\la)^2s}$ and so
by redefining \eqref{eq:char-mn} and \eqref{eq:zstar} with $s$ and $\la$, we have that the event $v^\la(\fl{\la^2s})<\fl{\la^2(1-\vr)^2\vr^{-2}s}$ is equivalent to
\[
\ba
Z^{*\la}(s)&=[\fl{(1-\la)^2s}-v^\la(\fl{\la^2s})]^+-[\fl{\la^2s}-w^\la(\fl{(1-\la)^2s})]^+\\
&=\fl{(1-\la)^2s}-v^\la(\fl{\la^2s})\\
&>\fl{(1-\la)^2s}-\fl{\la^2(1-\vr)^2\vr^{-2}s}.
\ea
\]
By $Z^{*\la}\overset{\text{d}}{=}Z^\la$, we conclude
\be
\begin{split}
&\Pv\bigl\{ v^\la(\fl{\vr^2t})<\fl{(1-\vr)^2t}\bigr\} \\
&\qquad
=\Pv\bigl\{Z^\la(s)>\fl{(1-\la)^2s}-\fl{\la^2(1-\vr)^2\vr^{-2}s}\bigr\}.
\end{split}
\label{eq:vtoZ}
\ee
Utilizing the definitions \eqref{eq:def-temp-la} and
 \eqref{eq:def-temp-r-1} of $\la$ and  $r$, one can
check that by  increasing $t_0=t_0(a,b,\vr)$ if necessary, one can guarantee that for
 $t\geq t_0$  there exists
a constant $C=C(\vr)$ such that
\[
\fl{(1-\la)^2s}-\fl{\la^2(1-\vr)^2\vr^{-2}s}\geq Crs^{2/3}.
\]
Combining this  with Corollary \ref{cor:Z-ub} and definition
\eqref{eq:def-temp-r-1} of $r$  we get the bound
\[
\begin{split}
\Pv\bigl\{ v^\la(\fl{\vr^2t}) <\fl{(1-\vr)^2t}\bigr\}
&\leq \Pv\bigl\{ Z^\la(s) > Crs^{2/3} \bigr\} \\
&\leq Cr^{-3} \leq C (\va/b)^3.
\end{split}
\]

Returning to \eqref{eq:lb-temp-2} to combine
all the bounds, we have
\[
\Pv\Bigl\{ \,\sup_{1\leq z\leq \va t^{2/3}}\bigl(  U^\vr_z + A_z(t) -A_1(t)\bigr)
\geq bt^{1/3} \Bigr\}
 \leq  C\Bigl( \,\frac{\va^3}{b^3}+ \frac{\va^3}{b^6}\,\Bigr).
\qedhere \]
\end{proof}

\begin{lm}  We have the asymptotics
\[
\lim_{\ve\searrow 0} \limsup_{t\to\infty}
\Pv\{ 0< U^\vr_{Z^{\vr}(t)^+} \leq \ve t^{2/3} \} =0.
\]
\label{lm:lb-2}
\end{lm}

Note that part of the event is the requirement $Z^\vr(t)>0$.

\begin{proof}
The limit comes from control over the point $Z^\vr(t)$. First
write
\begin{align*}
\Pv\{ 0< U^\vr_{Z^{\vr}(t)^+} \leq \ve t^{2/3} \}
\leq \Pv\{ 0<Z^\vr(t)\leq \delta t^{2/3}\}
+ \Pv\{ U^\vr_{\tfl{\delta t^{2/3}}}\leq \ve t^{2/3} \}.
\end{align*}
Given $\delta>0$, the last probability vanishes as $t\to\infty$
for any $\ve< \delta(1-\vr)^{-1}$. Thus it remains to show that
the first probability on the right can be made  arbitrarily small
for large $t$, by choosing a small enough $\delta$.

Let $0<\delta,b<1$.
\begin{align}
&\Pv\{ 0<Z^\vr(t)\leq \delta t^{2/3}\} \nn\\
&\leq \Pv\Bigl\{ \,\sup_{x>\delta t^{2/3}}\bigl( U^\vr_x+A_x(t)\bigr) <
\sup_{1\leq x\leq \delta t^{2/3}}\bigl( U^\vr_x+A_x(t)\bigr) \Bigr\}\nn\\
&\leq \Pv\Bigl\{ \,\sup_{x>\delta t^{2/3}}
\bigl( U^\vr_x+A_x(t)-A_1(t)\bigr) < bt^{1/3}\Bigr\}\label{eq:prob-temp-4}\\
&\qquad +\;
\Pv\Bigl\{ \,\sup_{1\leq x\leq \delta t^{2/3}}
\bigl( U^\vr_x+A_x(t)-A_1(t)\bigr) >  bt^{1/3}\Bigr\}.
\label{eq:prob-temp-5}
\end{align}

By Lemma \ref{lm:lb-1}
 the probability \eqref{eq:prob-temp-5} is bounded
by $C\delta^3(b^{-3}+b^{-6})$.
Bound the probability \eqref{eq:prob-temp-4} by
\begin{align}
&\Pv\Bigl\{ \,\sup_{\delta t^{2/3}<x\leq t^{2/3}}
\bigl( U^\vr_x+A_x(t)-A_1(t)\bigr) < bt^{1/3}\Bigr\}\nn\\
&\quad \leq
\Pv\bigl\{ v^\la(\fl{\vr^2t})> \fl{(1-\vr)^2t}-t^{2/3}\bigr\}
\label{eq:prob-temp-6}\\
&\qquad  + \Pv\Bigl\{ \,\sup_{\delta t^{2/3}<x\leq t^{2/3}}
(  U^\vr_x -V^\la_{x})
< bt^{1/3} \Bigr\} \label{eq:prob-temp-7}
\end{align}
where, following the example of
 the previous proof, we have introduced a new density,
this time
\[
\la=\vr+rt^{-1/3},
\]
and then used the reversal trick
of equation \eqref{eq:reversal-1} and Lemma \ref{lm:coup-3}
to deduce
\[
A_x(m,n)-A_1(m,n)\geq G^\la(m-x+1,n)-G^\la(m,n)\equiv -V^\la_{x-1}
\geq -V^\la_x,
\]
whenever $v^\la(\fl{\vr^2t})\leq \fl{(1-\vr)^2t}-t^{2/3}$. We claim that, given $\eta>0$ and parameter  $r$ from above,
we can fix  $\delta, b>0$  small
enough so that, for some $t_0<\infty$,  the probability in
\eqref{eq:prob-temp-7} satisfies
\be
\Pv\Bigl\{ \,\sup_{\delta t^{2/3}<x\leq t^{2/3}}
(  U^\vr_x -V^\la_{x})
< bt^{1/3} \Bigr\} \leq \eta
\qquad\text{for all $t\geq t_0$.}
\label{eq:lb-temp-11}
\ee
As $t\to\infty$,
\begin{align*}
&t^{-1/3} \Ev(U^\vr_{\fl{yt^{2/3}}} -V^\la_{\fl{yt^{2/3}}})
\longrightarrow \frac{-ry}{(1-\vr)^2}\\
&\quad\quad\text{and}\quad
t^{-2/3}\Vv(U^\vr_{\fl{yt^{2/3}}} -V^\la_{\fl{yt^{2/3}}})
\longrightarrow \frac{2y}{(1-\vr)^2}\equiv\sigma^2(\vr)y
\end{align*}
uniformly over $y\in[\delta, 1]$. Since we have a sum of
i.i.d's,
the probability in \eqref{eq:lb-temp-11} converges, as $t\to\infty$, to
\[
\Pv\Bigl\{ \,\sup_{\delta \leq y\leq 1}\Bigl(
\sigma(\vr) B(y)-\frac{ry}{(1-\vr)^{2}}\Bigr)
\leq  b \Bigr\}
\]
where $B(\cdot)$ is standard Brownian motion. The random variable
\[ \sup_{0\leq y\leq 1}\Bigl(
\sigma(\vr) B(y)-\frac{ry}{(1-\vr)^{2}}\Bigr)
\]
is positive almost surely, so the above probability is less
than $\eta/2$ for small $\delta$ and $b$. This implies
\eqref{eq:lb-temp-11}.

The probability in \eqref{eq:prob-temp-6}
is bounded by
\begin{align*}
&\Pv\bigl\{ v^\la(\tfl{\vr^2t})>\tfl{(1-\vr)^2t-t^{2/3}-1}\bigr\}\\
&\leq  \Pv\bigl\{ v^{1-\la}(\tfl{(1-\vr)^2t-t^{2/3}-1})<\tfl{\vr^2t}\bigr\}\\
&\leq \Pv\bigl\{ v^{1-\la}(\tfl{(1-\la)^2s})<\tfl{\la^2s}- qs^{2/3}\bigr\}\\
&=\Pv\bigl\{ Z^{*1-\la}(s)>  qs^{2/3}\bigr\}\\
&=\Pv\bigl\{ Z^{1-\la}(s)>  qs^{2/3}\bigr\}\\
&\leq Cq^{-3}.
\end{align*}
Above we first used  \eqref{eq:v-w-ineq} and
 transposition of the array $\{\om_{ij}\}$. Because this exchanges
the axes, density $\la$ becomes $1-\la$.
 Then we defined
$s$  by
\[
(1-\la)^2s= (1-\vr)^2t-t^{2/3}-1
\]
and observed that for large enough $t$, the second inequality
holds for some $q=
C(\vr)((1-\vr)r-\vr )$. We used \eqref{eq:zstar} and the distributional identity of $Z$ and $Z^*$ thereafter. The last inequality is
from Corollary \ref{cor:Z-ub}.

Now given $\eta>0$, choose $r$ large enough so that
$Cq^{-3}<\eta$.  Given this $r$, choose $\delta, b$ small enough so
that \eqref{eq:lb-temp-11} holds.  Finally, shrink $\delta$ further
so that $C\delta^3(b^{-3}+b^{-6})<\eta$ (shrinking $\delta$ does not
violate \eqref{eq:lb-temp-11}).
To summarize, we have shown that, given $\eta>0$,
if $\delta$ is small enough, then for all large  $t$
\be
\Pv\{1\leq  Z^\vr(t)\leq \delta t^{2/3}\}\leq 3\eta.
\label{eq:b-proof}
\ee
This concludes the proof of the lemma.
\end{proof}

Via transpositions we get the previous lemma also for the j-axis:

\begin{cor}   We have the asymptotics
\[
\lim_{\ve\searrow 0} \limsup_{t\to\infty}
\Pv\{ 0< U^\vr_{-Z^{\vr}(t)^-} \leq \ve t^{2/3} \} =0.
\]
\label{cor:lb-3}
\end{cor}
\begin{proof} Let  $\{\om_{ij}\}$ be an initial assignment
with density $\vr$. Let $\wt\om_{ij}=\om_{ji}$ be the
transposed array,
which is an initial assignment with density $1-\vr$.
Under transposition the
$\fl{(1-\vr)^2t}\times\fl{\vr^2 t}$ rectangle has become
$\fl{\vr^2 t}\times\fl{(1-\vr)^2t}$, the correct characteristic
dimensions for density $1-\vr$.  Since  transposition exchanges
the coordinate axes, after transposition
$U^\vr_{Z^\vr(t)^+}$ has become $U^{1-\vr}_{-Z^{1-\vr}(t)^-}$,
and so these two random variables
have the same distribution. The corollary is now
a consequence of Lemma \ref{lm:lb-2} because this lemma is valid
for each density $0<\vr<1$.
\end{proof}

\eqref{eq:b-proof} proves part (b) of Theorem \ref{tm:Z-1}.
The theorem below gives the lower bound for Theorem \ref{tm:G-var-1}
and thereby completes its proof.

\begin{tm}
\[
\liminf_{t\to\infty}\frac{\Ev(U^\vr_{Z^\vr(t)^+})}{t^{2/3}}>0,
\quad\text{and}\quad\liminf_{t\to\infty}\frac{\Vv(G^\vr(t))}{t^{2/3}}>0.
\]
\end{tm}
\begin{proof}
Suppose
there exists a density $\vr$ and a sequence $t_k\to\infty$ such that
$t_k^{-2/3}\Vv(G^\vr(t_k))\to 0$.  Then by Lemma \ref{lm:ns}
\[
\frac{\Ev(U^\vr_{Z^\vr(t_k)^+})}{t_k^{2/3}}\to 0
\quad\text{and}\quad
\frac{\Ev(U^\vr_{-Z^\vr(t_k)^-})}{t_k^{2/3}}\to 0.
\]
From this and Markov's inequality
\[
\Pv\{ U^\vr_{Z^\vr(t_k)^+} >\ve t_k^{2/3}\}\to 0
\quad\text{and}\quad
\Pv\{ U^\vr_{-Z^\vr(t_k)^-} >\ve t_k^{2/3}\}\to 0
\]
for every $\ve>0$. This together with Lemma \ref{lm:lb-2}
and Corollary \ref{cor:lb-3} implies
\[
\Pv\{ U^\vr_{Z^\vr(t_k)^+} > 0\}\to 0
\quad\text{and}\quad
\Pv\{ U^\vr_{-Z^\vr(t_k)^-} > 0\}\to 0.
\]
But these statements imply that
\[
\Pv\{ Z^\vr(t_k) > 0\}\to 0
\quad\text{and}\quad
\Pv\{ Z^\vr(t_k) < 0\}\to 0,
\]
which is  a contradiction since these two probabilities
add up to 1 for each fixed $t_k$.
This proves the second claim of the theorem.

The first claim follows because it is equivalent to the second.
\end{proof}

\section{Rarefaction boundary conditions}
\label{sc:gen-bd}
In this section we prove results on the longitudinal and transversal fluctuations of a maximal path under more general
boundary conditions.
Abbreviate as before
\[ (m,n) = \bigl( \fl{(1-\vr)^2t},\fl{\vr^2t}\bigr).\]
We start by studying $A_0(t)=A_0(m,n)$, the maximal path to $(m,n)$ when there are no weights on the axes. We still use the boundary conditions \eqref{eq:bondis}, so that we have coupled $A_0(t)$ and $G^\vr(t)$.
We prove another version
of  Lemma \ref{lm:lb-1} to make it applicable for all $t\ge 1$.
\begin{lm}\label{lem:A0}
Fix $0<\al <1$. There exists a constant $C=C(\al, \vr)$
 such that, for each $t\ge 1$ and $b\ge C$,
\[ \Pv\{G^\vr(t) - A_0(t)\geq bt^{1/3}\} \leq Cb^{-3\al/2}.\]
\end{lm}
\begin{proof}
Note that
\be
\ba
\Pv\{G^\vr(t) - A_0(t)\geq bt^{1/3}\}
\leq & \ \Pv\{\sup_{|z|\leq at^{2/3}} U^\vr_z(t) + A_z(t) - A_0(t)\geq bt^{1/3}\}\\
& + \Pv\{\sup_{|z|\leq at^{2/3}} U^\vr_z(t) + A_z(t)\neq G^\vr(t)\}.
\label{eq:splitAG}
\ea
\ee
The last term of \eqref{eq:splitAG} can easily be dealt with using
Corollary \ref{cor:Z-ub}: there exists a $C=C(\vr)$ such that
\be
\ba
&\Pv\{\sup_{|z|\leq at^{2/3}} U^\vr_z(t) + A_z(t)\neq G^\vr(t)\}
\leq  \Pv\{ Z^\vr(t)\geq at^{2/3}\} \\
&\qquad  + \Pv\{ Z^\vr(t)\leq -at^{2/3}\}  \leq Ca^{-3}. \label{eq:ubneq}
\ea
\ee
For the first term of \eqref{eq:splitAG} we will use the results from the proof of Lemma \ref{lm:lb-1}. We split the range of $z$ into $[1,at^{2/3}]$ and $[-1,-at^{2/3}]$ and consider for now only the first part. Define
\[ \la = \vr - rt^{-1/3}.\]
We can use \eqref{eq:lb-temp-2} and \eqref{eq:lb-temp-2-b},
where we choose $a=b^{\al/2}$, $p=2$,  and $r=b^{\al/2}$.
Choose $C=C(\al, \vr)>0$ large enough so that
for $b\ge C$
\eqref{eq:lb-temp-2-a} is satisfied and the denominator
of the last bound in \eqref{eq:lb-temp-2-b} is at least
$b/2$. Then we can claim that, for all $b\ge C$ and $t\ge 1$,
\be
\ba
&\Pv\{\sup_{1\leq z\leq at^{2/3}} U^\vr_z(t) + A_z(t) - A_0(t)
\geq bt^{1/3}\} \\
& \qquad \leq \Pv\bigl\{ v^\la(\fl{\vr^2t})<\fl{(1-\vr)^2t}\bigr\}
 + Cb^{\al/2-2}.
\ea
\label{eq:sec7-temp-5}
\ee
From \eqref{eq:vtoZ} we get with $s=(\vr/\la)^2t$
\[
\Pv\bigl\{ v^\la(\fl{\vr^2t})<\fl{(1-\vr)^2t}\bigr\}=\Pv\bigl\{Z^\la(s)>\fl{(1-\la)^2s}-\fl{\la^2(1-\vr)^2\vr^{-2}s}\bigr\}.
\]
Now we continue differently  than in Lemma
\ref{lm:lb-1} so that $t$ is not forced to be large.
An elementary calculation yields
\[
\fl{(1-\la)^2s}-\fl{\la^2(1-\vr)^2\vr^{-2}s}\geq
2\,\frac{1-\vr}{\vr}\,rs^{2/3} + \frac{2\vr-1}{\vr^2}\,r^2s^{1/3}-1.
\]
We want to write down conditions
under which  the right-hand side above is at least
$\delta rs^{2/3}$ for some constant $\delta$ and all
$s\ge 1$.  First increase the above constant $C=C(\al, \vr)$
so that if $b=r^{2/\al}\ge C$, then
\[
\frac{1-\vr}{\vr}\,rs^{2/3} -1 \ge \frac{1-\vr}{2\vr}\,rs^{2/3}
\qquad\text{ for all $s\ge 1$. }
\]
Then choose $\eta=\eta(\al,\vr)>0$ small enough such that
whenever $b\in[C,\,\eta t^{2/(3\al)}]$
(in this case $r$ is small enough compared to $t^{1/3}$, but notice that the interval might as well be empty when $t$ is small),
\[
\frac{1-\vr}{\vr}\,rs^{2/3} \ge -\,\frac{2\vr-1}{\vr^2}\,r^2s^{1/3}.
\]
This last condition is vacuously true if $\vr\ge 1/2$.

Now we have for $C\le b\leq \eta t^{2/(3\al)}$ and with
$\delta=(1-\vr)/(2\vr)$,
\[ \fl{(1-\la)^2s}-\fl{\la^2(1-\vr)^2\vr^{-2}s}\geq \delta rs^{2/3}
\qquad\text{ for all $s\ge 1$. }
\]
If we combine this with \eqref{eq:sec7-temp-5} and
 Corollary \ref{cor:Z-ub},
we can state that for all $C\le b\le \eta t^{2/(3\al)}$
and $t\ge 1$,
\[
\Pv\{\sup_{1\leq z\leq at^{2/3}} U^\vr_z(t) + A_z(t) - A_0(t)
\geq bt^{1/3}\}  \le Cb^{-3\al /2}
 + Cb^{\al/2-2}.
\]
Same argument works (or just apply transposition)
 for the values $-at^{2/3}\le z\le 1$,
so this same upper bound is valid for the first probability
on the right-hand side of \eqref{eq:splitAG}.

Taking \eqref{eq:ubneq} also into consideration,
 at this point we have shown
that whenever $C\le b\le \eta t^{2/(3\al)}$ and $t\ge 1$,
\[ \Pv\{G^\vr(t) - A_0(t)\geq bt^{1/3}\} \leq \frac12C(b^{\al/2-2} +
b^{-3\al /2})\leq Cb^{-3\al /2}.\]

What if $b\geq \eta t^{2/(3\al)}$? Note that
\[ \Pv\{G^\vr(t) - A_0(t)\geq bt^{1/3}\} \leq \Pv\{G^\vr(t) \geq
bt^{1/3}\}.\]
Since $G^\vr(t)$ is the sum of two (dependent) random variables,
each of which in turn is the sum of i.i.d. exponentials, and since
\[
\Ev(G^\vr(t)b^{-1}t^{-1/3}) \leq C(\vr, \eta )b^{\al -1}
\]
($\Ev(G^\vr(t))$ is basically linear in $t$ by \eqref{eq:char-mn} and Lemma \ref{lm:NE}),
we conclude that $\Pv\{G^\vr(t) - A_0(t)\geq bt^{1/3}\}$ goes to
zero faster than any polynomial in $b$, if $b\geq \eta t^{2/(3\al)}$.
This proves the lemma for all $b\ge C$.
\end{proof}

Now we can establish that the
 fluctuations of $A_0(t)$ are of order $t^{1/3}$.

\begin{cor}\label{cor:At}
Fix $0<\al <1$. There exists a constant $C=C(\al, \vr)$ such
that for all $a>0$ and $t\ge 1$,
\[ \Pv \{ |A_0(t)-t| > at^{1/3}\} \leq Ca^{-3\al/2}.\]
In particular this means that
\[ \Ev(|A_0(t) - t|) = O(t^{1/3})\quad\text{and}\quad \Ev(A_0(t))=t -
O(t^{1/3}).\]
\end{cor}
\begin{proof}
Lemma \ref{lem:A0} together with Theorem \ref{tm:EUbd} implies for
$a\ge C(\al,\vr)$
\[
\ba
\Pv \{ |A_0(t)-t| > at^{1/3}\} &\leq \Pv\{G^\vr(t)-A_0(t)> at^{1/3}/2\}\\
& \ \ \  + \Pv \{ |G^\vr(t)-t| > at^{1/3}/2\}\\
& \leq C_1a^{-3\al/2} + C_2a^{-2}\\
& \leq Ca^{-3\al/2}.
\ea
\]
Finally, we can always increase $C$ in order to take all $0<a\leq
C(\al,\vr)$ values into account.
\end{proof}

We can also consider the fluctuations of the position of a maximal path.
To this end we extend the definition of $Z(t)$, the exit point from the
axes. We define $Z_l(t)$ as the $i$-coordinate of the right-most point on
the horizontal line $j=l$ of the right-most maximal path to $(m,n)$ (we
say right-most path, because later in this section we will consider
boundary conditions that no longer necessarily have a unique longest
path). We will use the notation $Z^\vr_l$ to denote the stationary
situation and $Z^0_l$ to denote the situation where all the weights
on the axes are zero. Note that in all cases
\[ Z(t)^+ = Z_0(t).\]
\begin{lm}
Define $(k,l)=(\fl{(1-\vr)^2s},\fl{\vr^2s})$ for $s\leq t$. There
exists a constant $C=C(\rho)$ such that for all $s\leq t$ with
$t-s\geq 1$ and all $a>0$
\[ \Pv\{ Z^\vr_l(t)\geq k+a(t-s)^{2/3}\} \leq Ca^{-3}.\]
\end{lm}
\begin{proof}
There are several ways to see this, for example using time-reversal.
One can also pick a new origin at $(k,l)$, and define a
  last-passage model in
the rectangle $[k,m]\times[l,n]$ with  boundary conditions given
by $I$- and $J$-increments of the  $G$-process in the original
rectangle $[0,m]\times[0,n]$.  The maximizing path in this
new model connects up with the original maximizing path.
  Hence in this new model
 it looks as though the maximal path to
$(m-k,n-l)$ exits the $i$-axis beyond the point $a(t-s)^{2/3}$, and so
\[ \Pv\{ Z^\vr_l(t)\geq k+a(t-s)^{2/3}\} = \Pv\{ Z^\vr(t-s)\geq a(t-s)^{2/3}\}.\]
We have ignored the integer parts here, but this can be dealt with uniformly in $a>0$. Now we can use Corollary \ref{cor:Z-ub} to conclude that
\[ \Pv\{ Z^\vr_l(t)\geq k+a(t-s)^{2/3}\} \leq Ca^{-3}.\]
\end{proof}

To get a similar result for $Z^0_l(t)$ we need a more convoluted
argument and the conclusion is  a little weaker.

\begin{lm}\label{lem:Z0sm}
Define $(k,l)=(\fl{(1-\vr)^2s},\fl{\vr^2s})$ for $s\leq t$. There exists a constant $C=C(\al, \vr)$ such that for all $a>0$ and $t\geq 1$
\[ \Pv\{ Z^0_l(t)> k+a\,t^{2/3}\} \leq Ca^{-3\al }.\]
\end{lm}
\begin{proof}
The event $\{Z^0_l(t)> k+u\}$ is equivalent to the event
\[ E=\left\{A_0(t) = \sup_{z\geq u}\{ A_0(k+z+1,l) + \tilde{A}_0(z-u,0)\}\right\}.\]
Here, $A_0(i,j)$ is the weight of the maximal path (not using the axes) from $(0,0)$ to $(i,j)$, including the endpoint, whereas $\tilde{A}_0(i,j)$ is the weight of the maximal path from $(k+u+i,l+j)$ to $(m,n)$, including the endpoint but excluding the starting point {\em and excluding all the weights directly to the right or directly above $(k+u+i,l+j)$}. This corresponds to choosing $(k+u+i,l+j)$ as a new origin, and making sure that the axes through this origin have no weights. Note that the processes $A_0(\cdot,l)$ and $\tilde{A}_0(\cdot,0)$ are independent. The idea is to
bound $A_0$ and $\tilde{A}_0$ by appropriate stationary processes
 $G^\la$ and $G^{\tilde{\la}}$ to show that,
 with high probability, this supremum will be too small if $u$ is too
large.   We can couple the processes $G^\la$ and $G^{\tilde{\la}}$, where $\tilde{\la}>\la$, in the following way: $G^\la$ induces weights on the horizontal line $j=l$ through the increments of $G^\la$, see Lemma \ref{lm:NE}. The process $G^{\tilde{\la}}$ takes the point
$(k+u,l)$ as origin and uses as boundary weights
on the horizontal line $j=l$, with a slight abuse of
notation, for $i\geq 1$
\[ \tilde{\om}_{i0} = \frac{1-\la}{1-\tilde{\la}} I_{k+u+i+1,l}= \frac{1-\la}{1-\tilde{\la}}(G^\la(k+u+i+1,l)-G^\la(k+u+i,l)).\]
These weights are independent $\Expd(1-\tilde{\la})$ random variables. The weights $\tilde{\om}_{0j}\sim
\Expd(\tilde{\la})$ on the line $i=k+u$ can be chosen
independently of everything else, whereas for $i,j\geq 1$
\[ \tilde{\om}_{ij} = \om_{u+k+i,l+j}.\]
So $G^{\tilde{\la}}(i,j)$ equals the weight of the maximal path
from $(k+u,l)$ to $(k+u+i,l+j)$, using as weights on the
points $(k+u+i,l)$ the $\tilde{\om}_{i0}$ (for $i\geq 1$), on the
points $(k+u,l+j)$ the $\tilde{\om}_{0j}$ (for $j\geq 1$) and on
the points $(k+u+i,l+j)$ the original $\om_{k+u+i,l+j}$
(for $i,j\geq 1$). This construction leads to
\[ A_0(i,j)\leq G^\la(i,j)\quad \text{and}\quad
\tilde{A}_0(i,j)\leq G^{\tilde{\la}}(m-k-u,n-l) - G^{\tilde{\la}}(i,j).\]
Also, for all $0\leq i\leq m-k-u-1$,
\[ G^\la(k+u+i+1,l)-G^\la(k+u+1,l)\leq G^{\tilde{\la}}(i,0).\]
Therefore, for all $0\leq i\leq m-k-u-1$,
\[
\ba
A_0(k+u+i+1,l) + \tilde{A}_0(i,0) &\leq G^\la(k+u+i+1,l)- G^{\tilde{\la}}(i,0) + \\
&\ \ \ \ G^{\tilde{\la}}(m-k-u,n-l) \\
&\leq G^\la(k+u+1,l) + G^{\tilde{\la}}(m-k-u,n-l).
\ea
\]
So we get
\be
\label{eq:boundPE}
\Pv (E) \leq \Pv\{ A_0(t) - G^\la(k+u+1,l) -
G^{\tilde{\la}}(m-k-u,n-l)\leq 0\}.
\ee
Here, we can still choose $\la$ and $\tilde{\la}$ as long as $0<\la <\tilde{\la}$, but it is not hard to see that for the optimal choices (in expectation) of $\la$ and $\tilde\la$ are determined by
\be
\frac{(1-\la)^2}{\la^2}=\frac{k+u+1}{l}\quad\text{and}\quad\frac{(1-\tilde\la)^2}{\tilde\la^2}=\frac{m-k-u}{n-l}.\label{eq:optchar}
\ee
With these choices we get
\[ \Ev(G^\la(k+u+1,l)) = (\sqrt{k+u+1} + \sqrt{l})^2\]
and
\[ \Ev(G^{\tilde{\la}}(m-k-u,n-l)) = (\sqrt{m-k-u}+\sqrt{n-l})^2.\]
This particular choice of $(\la, \tilde{\la})$ is valid (i.e., $\tilde{\la}>\la$) as soon as $u\geq C(\vr)$. Smaller $u$ can be dealt with by increasing $C$ in the statement of lemma. We have for $u\geq 2$
\[ 
\ba
\Ev(G^\la(k+&u+1,l) + G^{\tilde{\la}}(m-k-u,n-l)) =  m+n+2\sqrt{l(k+u+1)}\\
&  \hspace{5cm} +2\sqrt{(n-l)(m-k-u)}+1\\
& \leq  ((1-\vr)^2+\vr^2)t + 1 + 2\sqrt{\vr^2s}\,\sqrt{(1-\vr)^2s+u} + \sqrt{\frac{l}{k+u}}\\
& \hspace{2.5cm} + 2\sqrt{\vr^2(t-s)+1}\,\sqrt{(1-\vr)^2(t-s)-u+1}\\
& \leq  ((1-\vr)^2+\vr^2)t + C(\vr) + 2\sqrt{\vr^2s}\,\sqrt{(1-\vr)^2s+u}\\
& \hspace{2cm} + 2\sqrt{\vr^2(t-s)}\,\sqrt{(1-\vr)^2(t-s)-(u-1)}\\
& \leq t + C(\vr) + \frac{\vr}{1-\vr}u - \frac{\vr}{1-\vr}(u-1) - \frac14\, \frac{\vr}{(1-\vr)^3}\, \frac{(u-1)^2}{t-s}\\
& \leq t - C_1(\vr)\,\frac{u^2}{t} + C_2(\vr).
\ea
\]
If $u=\fl{at^{2/3}}$, then we can choose constants $M=M(\vr)$ and $C_1=C(\vr)$ such that for all $a>M$ and $t\geq 1$,
\be \label{eq:EGla}
\Ev(G^\la(k+u+1,l) + G^{\tilde{\la}}(m-k-u,n-l)) \leq t - C_1 a^2t^{1/3}.
\ee
Smaller $a$ can be dealt with by increasing the constant $C$ in the statement of the lemma. Now note that, using (\ref{eq:boundPE}), we get 
\[
\ba
\Pv(E) &\leq  \Pv\{ A_0(t) - t\leq G^\la(k+u+1,l) + G^{\tilde{\la}}(m-k-u,n-l) - t\}\\
&\leq  \Pv( A_0(t)-t \leq -\frac12C_1a^2t^{1/3}) \\
& \ \ \    + \Pv( G^\la(k+u+1,l) + G^{\tilde{\la}}(m-k-u,n-l) - t\geq -\frac12 C_1a^2t^{1/3})\\
&\leq  \Pv( A_0(t)-t \leq -\frac12 C_1a^2t^{1/3}) + C_2a^{-4} \leq Ca^{-3\al}.
\ea
\]
For the last line we used (\ref{eq:EGla}), the fact that
\[ \Vv(G^\la(k+u+1,l) + G^{\tilde{\la}}(m-k-u,n-l))\leq Ct^{2/3},\]
(notice that the choice \eqref{eq:optchar} places these coordinates in the $G$'s on the respective characteristics, see \eqref{eq:char-mn}), and Corollary \ref{cor:At}.

\end{proof}

We now turn to the case of a rarefaction fan introduced by \eqref{eq:rareom}.
\begin{proof}[Proof of Theorem \ref{tm:rare}]
The statement follows from the trivial observation that
\[ A_0(t)\leq \hat{G}(t) \leq G^{\rho}(t)\]
(if there is less weight, the paths get shorter), Corollary \ref{cor:At} and Theorem \ref{tm:EUbd}.
\end{proof}
\begin{proof}[Proof of Theorem \ref{tm:raretrans}]
For the first inequality, we introduce the process $G^{\Wc=0}$, which uses the same weights as $G^\rho$, except on the $j$-axis, where all weights are 0 (so $\om^{\Wc=0}_{0j}=0$). It is not hard to see that
\[ \hat{Z}_l(t) \leq Z^{\Wc=0}_l(t),\]
simply because the right-most maximal path for $G^{\Wc=0}$ stays at least as long on the $i$-axis as a maximal path for $\hat{G}$, and it can coalesce with, but never cross a maximal path for $\hat{G}$. So we get
\[ \Pv\{ \hat{Z}_l(t)\geq k+at^{2/3}\} \leq \Pv\{ Z^{\Wc=0}_l(t)\geq k+at^{2/3}\}.\]
First we will show that with high probability, $Z^{\Wc=0}_0(t)$ is not too large. This will imply that if $Z^{\Wc=0}_l(t)$ is large, it must be because $Z^0_l(\tilde{t})$ (for an appropriately chosen $\tilde{t}$) is large, which has low probability because of Lemma \ref{lem:Z0sm}.

Note that, as in the proof of the previous lemma,
\be\label{eq:ZW>u}
\{ Z^{\Wc=0}_0(t)> u\} = \{ G^{\Wc=0}(t) = \sup_{z>u} (U^\vr_z + A_z(t))\}.
\ee
Now define a stationary process $G^\la$, with $\la > \vr$, whose origin is placed at $(u,0)$. It uses as weights on the $i$-axis
\[ \om^\la_{i0} = \frac{1-\vr}{1-\la}\om_{u+i+1,0}.\]
On the line $i=u$, $G^\la$ uses independent $\Expd({\la})$ weights. This construction guarantees that for $i\geq 0$
\[ U^\vr_{u+i+1} - U^\vr_{u+1} \leq G^\la(i,0).\]
Also, for $z>u$,
\[ A_z(t) \leq G^\la(m-u,n)-G^\la(z-u-1,0).\]
This implies that
\[ \sup_{z>u} (U^\vr_z + A_z(t))\leq U^\vr_{u+1} + G^\la(m-u,n).\]
This means that, using (\ref{eq:ZW>u}),
\be\label{eq:ZWsubset}
\{ Z^{\Wc=0}_0(t)> u\} \subset \{ G^{\Wc=0}(t) \leq  U^\vr_{u+1} + G^\la(m-u,n)\}.
\ee
Again we have that for the optimal $\la$,
\[ \Ev(G^\la(m-u,n)) = (\sqrt{m-u}+\sqrt{n})^2,\]
which leads to
\[ 
\ba
\Ev(U^\vr_{u+1} &+ G^\la(m-u,n)) \leq \frac{u+1}{1-\vr} + m + n - u + 2\sqrt{n}\,\sqrt{m-u}\\
& \leq (1-\vr)^2t + \vr^2 t +\frac{\vr u}{1-\vr} + C_1(\vr) + 2\sqrt{\vr^2t}\,\sqrt{(1-\vr)^2t-(u-1)}\\
& \leq t + C_2(\vr) -\frac14\, \frac{\vr}{(1-\vr)^3}\,\frac{(u-1)^2}{t}.
\ea
\]
Just as in the proof Lemma \ref{lem:Z0sm}, we see that if $u=\fl{bt^{2/3}}$, we can choose constants $M=M(\vr)$ and $C_1=C_1(\vr)$ such that for all $b>M$ and $t\geq 1$,
\[ \Ev(U^\vr_{u+1} + G^\la(m-u,n)) \leq t - C_1b^2t^{1/3}.\]
Note that with (\ref{eq:ZWsubset})
\[
\ba
\Pv( Z^{\Wc=0}_0(t)> bt^{2/3}) & \leq \Pv (G^{\Wc=0}(t) -t \leq -\frac12 C_1 b^2t^{1/3})\\
& \ \ \ + \Pv ( U^\vr_{u+1} + G^\la(m-u,n) -t \geq -\frac12 C_1 b^2t^{1/3}).
\ea
\]
Now we can use the fact that
\[ \Vv(U^\vr_{u+1} + G^\la(m-u,n)) = O(u+t^{2/3})\]
(again, the optimal choice for $\la$ has placed the coordinates in $G$ on the characteristics w.r.t.\ $\la$),
and Theorem \ref{tm:rare} to conclude that for $b>M$
\[ \Pv \{Z^{\Wc=0}_0(t)> bt^{2/3}\}\leq Cb^{-3\al}.\]
For $b\leq M$ we can increase $C$.

A little picture reveals that if $Z^{\Wc=0}_l(t)\geq k+at^{2/3}$ and $Z^{\Wc=0}_0(t)\leq at^{2/3}/2$, then the maximal path that does not use the weights on the axes from the point $(at^{2/3}/2,0)$ to $(m,n)$, must pass to the right of $(k+at^{2/3},l)$, an event with smaller probability than the event $\{ Z^0_l(t)\geq k + at^{2/3}/2\}$, which with Lemma \ref{lem:Z0sm} proves the first inequality.

The second inequality of the Theorem is a corollary of the first. We assume $k-at^{2/3}\ge0$, otherwise this statement is trivial. Also, we prove for $a>2\frac{(1-\vr)^2}{\vr^2}$, one can always increase $C$ if this is not the case. Fix $\wt s$ such that
\[
k-at^{2/3}=(1-\vr)^2\wt s,\quad\text{then}\quad k':\,=\fl{(1-\vr)^2\wt s},\qquad l':\,=\fl{\vr^2\wt s}.
\]
With these definitions,
\[
l\ge l'+\vr^2s-1-\frac{\vr^2}{(1-\vr)^2}\cdot(1-\vr)^2\wt s\ge l'+\frac{\vr^2}{(1-\vr)^2}\cdot at^{2/3}-1.
\]
Define also $\hat{Y}^T_{k'}$ to be the highest point of the left-most maximal path on the vertical line $i=k'$. As the left-most maximal path is North-East, we have
\[
\Pv\{\hat{Y}_l(t)\le k-at^{2/3}\}=\Pv\{\hat{Y}^T_{k'}\ge l\}\le\Pv\{\hat{Y}^T_{k'}\ge l'+\frac{\vr^2}{(1-\vr)^2}\cdot at^{2/3}-1\}.
\]
Pick $\wt a=a\vr^2/(1-\vr)^2-1>1$, then the right hand-side is bounded by $\Pv\{\hat{Y}^T_{k'}\ge l'+\wt at^{2/3}\}$. The transposed array $\wt\om_{ij}:\,=\om_{ji}$, $i,\,j\ge0$ has rarefaction fan boundary conditions w.r.t.\ the parameter $1-\vr$. Moreover, $\hat{Y}^T_{k'}$ becomes the right-most point of the right-most maximal path on the horizontal line $i=k'$ in the transpose picture. The first part of the Theorem with $1-\vr$, $\wt s\le t$ and $\wt a$ then completes the proof by $\wt a^{-3\al}<[\frac12(\wt a+1)]^{-3\al}=C'(\vr,\,\al)\cdot a^{-3\al}$.
\end{proof}
\bibliography{refsmarton}
\bibliographystyle{plain}

\bigskip
{\sc M.\ Bal\'azs, Mathematics Department, University of Wisconsin-Madison,} Van Vleck Hall, 480 Lincoln Dr, Madison WI 53706-1388, USA.\\
\indent
{\it E-mail address:} {\tt balazs@math.wisc.edu}

\bigskip
{\sc E. Cator, Delft University of Technology, faculty EWI,} Mekelweg 4, 2628CD, Delft, The Netherlands.\\
\indent
{\it E-mail address:} {\tt e.a.cator@ewi.tudelft.nl}

\bigskip
{\sc T.\ Sepp\"al\"ainen, Mathematics Department, University of Wis\-con\-sin-Madison,} Van Vleck Hall, 480 Lincoln Dr, Madison WI 53706-1388, USA.\\
\indent
{\it E-mail address:} {\tt seppalai@math.wisc.edu}

\end{document}